\documentclass[reqno]{amsart}
\usepackage{amsmath}
\usepackage{amssymb}

\theoremstyle{plain}
\newtheorem{theorem}{Theorem}

\newtheorem{lemma}{Lemma}

\theoremstyle{definition}

\theoremstyle{remark}

\numberwithin{equation}{section}

\input{tcilatex}

\begin{document}
\title[CALIBRATED DATA]{{\large \textbf{ANOMALIES IN THE ANALYSIS\\
OF CALIBRATED DATA}}}
\author{D. R. Jensen and D. E. Ramirez}
\address{Department of Mathematics\\
University of Virginia\\
P.O. Box 400137\\
Charlottesville, VA 22904--4137}
\email{der@virginia.edu}
\thanks{}
\date{}

\begin{abstract}
This study examines effects of calibration errors on model assumptions and
data--analytic tools in direct calibration assays. These effects encompass
induced dependencies, inflated variances, and heteroscedasticity among the
calibrated measurements, whose distributions arise as mixtures. These
anomalies adversely affect conventional inferences, to include the
inconsistency of sample means; the underestimation of measurement variance;
and the distributions of sample means, sample variances, and Student's 
\textit{t} as mixtures. Inferences in comparative experiments remain largely
intact, although error mean squares continue to underestimate the
measurement variances. These anomalies are masked in practice, as
conventional diagnostics cannot discern irregularities induced through
calibration. Case studies illustrate the principal issues.
\end{abstract}

\maketitle

{}

\section{Introduction}

\label{S:intro}

Calibrated measurements, intrinsic to the sciences and engineering, are
inherently subject to errors of calibration. These errors induce
dependencies in violation of a basic tenet in much of applied statistics,
namely, that observations should be uncorrelated if not independent. These
issues traditionally have been overlooked by both scientists and
statisticians, despite a century of emerging methodologies for the analysis
of experimental data. Not only are many parametric and nonparametric
procedures at risk under such violations, but so also are conventional
diagnostics for checking critical features of a model. We return to these
subsequently.

To fix ideas, observed responses $\{{Z}_{0} ,{Z}_{1},\ldots,{Z}_{n} \}$ are
often adjusted to ${Z}_{0} $ as standard, giving differences $\{{Y}_{i}  = ({%
Z}_{i}  - {Z}_{0} ); 1 \leq i \leq n\}$ as the objects of interest to the
investigator. Moreover, if $\{{Z}_{0} ,{Z}_{1},\ldots,{Z}_{n} \}$ are
mutually uncorrelated having variances $\{\QTR{up}{Var}({Z}_{0} ) = {\sigma }%
_{0}^{2} ,\,\QTR{up}{Var}({Z}_{i} ) = {\sigma }^{2} ; 1 \leq i \leq n\},$
then $\{{Y}_{1},\ldots,{Y}_{n} \}$ are equicorrelated with parameter $\rho = 
{\sigma }_{0}^{2} /({\sigma }^{2}  + {\sigma }_{0}^{2} ),$ having variances $%
\{\QTR{up}{Var}({Y}_{i} ) = {\sigma }^{2}  + {\sigma }_{0}^{2} ; 1 \leq i
\leq n\}$ that are inflated in comparison with unadjusted values.

Linearly calibrated instruments are pervasive. Some unintended consequences,
to be examined here, include (i) the structure of induced dependencies,
heteroscedasticity, and other departures from conventional model
assumptions; (ii) the inflation of measurement variances in comparison with
intended values; and (iii) effects of calibration on conventional inferences
for location, scale, and model diagnostics. We first examine moments, then
effects of calibration on the actual measurement distributions themselves.
We focus here on \textit{direct calibration assays} to be identified
subsequently. An outline follows.

Section 2 gives notation and technical support. Section 3 reconsiders the
calibration process with reference to irregularities induced through
calibration errors. Section 4 addresses the impact of these irregularities
on conventional inferences, to include (i) inferences regarding the mean and
variance in a single sample; and (ii) the analysis of one--way experimental
data, including tests on means and variances. The latter remain largely
intact, although measurement variances continue to be underestimated.
Section 5 reexamines the ability of conventional diagnostics to uncover
violations induced through calibration. Section 6 undertakes a case study to
illustrate essential findings. Section 7 offers a brief summary and
cautionary note. Some peripheral matters are deferred to an Appendix.

\section{Preliminaries}

\label{S:prelm}

\subsection{Notation}

\label{SS:nota} Designate $\mathbb{R}^{n}$ as Euclidean \textit{n}--space, $%
\mathbb{R}^{n}_{+}$ as its positive orthant, $\mathbb{S}_{n}$ as the real
symmetric $(n \times n)$ matrices, and $\mathbb{S}_{n}^{+}$ and $\mathbb{S}%
_{n}^{0}$ as their positive definite and positive semidefinite varieties.
Arrays appear in bold type, to include the transpose $\boldsymbol{A}
^{\prime}$ and inverse $\boldsymbol{A}^{-1} $ of $\boldsymbol{A} ;$ the unit
vector $\boldsymbol{1}_{n}  = [1,\ldots,1]^{\prime}\in\mathbb{R}^{n} ;$ the
identity matrix $\boldsymbol{I}_{n} ;$ a block--diagonal matrix $\QTR{up}{%
Diag}({\boldsymbol{A} }_{1},\ldots,{\boldsymbol{A} }_{k} );$ and $%
\boldsymbol{B}_{n} $ = $(\boldsymbol{I}_{n}  - n^{-1}\boldsymbol{1}_{n} 
\boldsymbol{1}_{n} ^{\prime}).$ Following Loewner (1934), matrices $(%
\boldsymbol{A} ,\boldsymbol{B} )$ in $\mathbb{S}_{n}$ are said to be ordered
as $\boldsymbol{A} \succeq_{L}\!\boldsymbol{B} $ for $\boldsymbol{A} -%
\boldsymbol{B} \in\mathbb{S}_{n}^{0},$ with $\boldsymbol{A} \succ_{L}\!%
\boldsymbol{B} $ whenever $\boldsymbol{A} -\boldsymbol{B} \in\mathbb{S}%
_{n}^{+}.$ Moreover, $\mathcal{C}(n)$ comprises the convex sets in $\mathbb{R%
}^{n}$ symmetric under reflection through $\boldsymbol{0} \in\mathbb{R}^{n}.$
Operators $E(\boldsymbol{Y} )$ and $\QTR{up}{V}(\boldsymbol{Y} )$ designate
the expected vector and dispersion matrix for $\boldsymbol{Y} \in\mathbb{R}%
^{n},$ with $\QTR{up}{Var}(Y)$ as the variance on $\mathbb{R}^{1}.$ We
further require $\{{\mu}_{r} (Z) = E(Z^{r}); r=1,2\}$ as moments about $0\in%
\mathbb{R}^{1},$ identifying ${\kappa }_{2}  = {\mu}_{2} ({\widehat{\beta}}%
_{1} )$ in terms of a linear estimatot ${\widehat{\beta}}_{1} $ to be
encountered subsequently. The comparative concentration of probability
measures on $\mathbb{R}^{n}$ may be gauged on defining the measure $%
\mu(\cdot)$ to be \textit{more peaked} about $\boldsymbol{0} \in\mathbb{R}%
^{n}$ than $\nu(\cdot),$ if and only if $\mu(A)\ge\nu(A)$ for every set $A\in%
\mathcal{C}(n),$ as in Sherman (1955). Specifically, the peakedness ordering
for scale mixtures of Gaussian measures on $\mathbb{R}^{n}$ is tantamount to
the stochastic ordering of their mixing distributions, as demonstrated in
the Appendix.

\subsection{Special Distributions}

\label{SS:dsns} Here \textit{pdf} and \textit{cdf} refer to probability
density and cumulative distribution functions; for $\boldsymbol{Y} \in%
\mathbb{R}^{n},$ $\mathcal{L}(\boldsymbol{Y} )$ designates its law of
distribution and $G(\boldsymbol{y}  )$ its \textit{cdf;} and \textit{iid}
refers to independent and identically distributed random elements.
Distributions of note on $\mathbb{R}^{1}$ include the Gaussian law $%
N_{1}(\mu,{\sigma }^{2} ),$ with parameters $(\mu,{\sigma }^{2} );$
noncentral versions of Student's $t(\nu,\lambda )$ and ${t}^{2} (\nu,\lambda
),$ chi--squared ${\chi}^{2} (\nu,\lambda ),$ and Snedecor--Fisher $F({\nu}%
_{1} ,{\nu}_{2} ,\lambda )$ distributions, with $\{\nu,{\nu}_{1} ,{\nu}_{2}
\}$ as degrees of freedom and $\lambda $ as a noncentrality parameter; and ${%
G}_{0} (\alpha ,\beta )$ as the gamma distribution on $\mathbb{R}^{1}_{+}$
having parameters $(\alpha ,\beta ).$ In particular, ${g}_{T} (t;\nu,\lambda
),$ ${g}_{T^{2}} (u;\nu,\lambda ),$ ${g}_{F} (u;{\nu}_{1} ,{\nu}_{2}
,\lambda ),$ and ${g}_{0} (u;\alpha ,\beta )$ designate the densities
corresponding to $t(\nu,\lambda ),$ ${t}^{2} (\nu,\lambda ),$ $F({\nu}_{1} ,{%
\nu}_{2} ,\lambda ),$ and ${G}_{0} (\alpha ,\beta),$ respectively.

To continue, $N_{n}(\boldsymbol{\theta} ,\boldsymbol{\Sigma } )$ designates
the Gaussian law on $\mathbb{R}^{n},$ and ${g}_{n} (\boldsymbol{y} ;%
\boldsymbol{\theta} ,\boldsymbol{\Sigma } )$ its \textit{pdf,} having
location--scale parameters $(\boldsymbol{\theta} ,\boldsymbol{\Sigma } ).$
Ensembles on $\mathbb{R}^{n},$ and mixtures over these, include the
translation--scale mixtures 
\begin{equation}
{f}_{1} (\boldsymbol{y} ;\boldsymbol{\theta} ,\boldsymbol{\Sigma } ,{G}_{1}
) = \int_{-\infty}^{\infty} {g}_{n} (\boldsymbol{y} ;\boldsymbol{\theta} (t),%
\boldsymbol{\Sigma } (t))d{G}_{1} (t),
\end{equation}
\noindent and purely scale mixtures when $\boldsymbol{\theta}  = \boldsymbol{%
0} .$ Nonstandard distributions for quadratic forms proceed conditionally on
letting $\mathcal{L}(U\!\mid\!w)$ have the scaled gamma density ${g}_{0}
(u;\alpha ,w\beta )$ = $(w\beta )^{-\alpha } u^{\alpha -1}e^{-u/w\beta
}/\Gamma(\alpha ),$ then compounding these as 
\begin{equation}
{f}_{3} (u;\alpha ,\beta ,{G}_{3} ) = \frac{u^{\alpha -1}}{\beta ^{\alpha
}\Gamma(\alpha )} \int_{0}^{\infty} w^{-\alpha } e^{-x/w\beta }d{G}_{3} (w)
\end{equation}
\noindent with ${G}_{3} (w)$ as a \textit{cdf} on $\mathbb{R}^{1}_{+}.$

Subsequent developments have links to exchangeable sequences. Consider $\{%
\boldsymbol{Z}_{0} , {Z}_{1} , {Z}_{2} , \ldots\}$ such that $\boldsymbol{Z}%
_{0} \in\mathbb{R}^{k}$ is independent of $\{{Z}_{1} , {Z}_{2} , \ldots\},$
whereas $\{{Z}_{1} , {Z}_{2} , \ldots\}$ are \textit{iid} on $\mathbb{R}%
^{\infty}.$ Further let $\psi : \mathbb{R}^{k+1}\!\rightarrow\!\mathbb{R}%
^{1};$ define $\{{Y}_{i}  = \psi(\boldsymbol{Z}_{0} ,{Z}_{i} ); i=1,
2,\ldots\};$ recall from DeFinetti's theorem that the sequence $\{{Y}_{1} ,{Y%
}_{2} ,\ldots\}$ is now \textit{exchangeable} on $\mathbb{R}^{\infty};$ and
infer that joint distributions projected onto $\mathbb{R}^{n}$ are invariant
under permutations. In short, $\{{Y}_{1},\ldots,{Y}_{n} \}$ are \textit{iid,}
under second moments having common values for the parameters $(\mu,{\sigma }%
^{2} ,\rho).$

\section{Calibration}

\label{S:calib}

\subsection{Essentials}

\label{SS:back}

Instruments are calibrated using outputs at successive inputs to determine a 
\textit{calibration curve;} new readings are assigned values on the scale of
measurements using the calibrated device; and periodic checks against a
standard determine when recalibration is required. In this study we utilize 
\textit{direct assays} in which instrumental readings $\{{X}_{i} ;1 \leq i
\leq n_{0}\}$ during calibration relate to observed measurements $\{{U}_{i}
;1 \leq i \leq n_{0}\}$ through $\{{U}_{i}  = {\beta}_{0}  + {\beta}_{1} {X}%
_{i}  + {\varepsilon }_{i} ; 1 \leq i \leq n_{0}\}.$ For example, the octane
rating $(U)$ in the production of gasoline relates linearly to the percent
of purity $(X)$ in a specimen to be assayed. Octane numbers require
expensive and time--consuming dynamic laboratory testing, whereas the
percent purity is readily determined. Once calibrated, the octane number of
a given specimen is determined vicariously from its percent purity. On the
other hand, \textit{indirect assays} proceed on reversing the roles of ${U}%
_{i} $ and ${X}_{i} $ during calibration. Models for calibration and their
analyses have been debated by several authors; for a summary and early
references see Krutchkoff (1971). Problems with moments and consistency
remain to be resolved in indirect assays, but the technical issues between
the two types differ mainly in detail. It is noteworthy that research has
yet to address the principal issues undertaken here, namely, irregularities
in models and supporting analyses attributable to calibration.

\subsection{Error Analysis}

\label{SS:erranal} To continue, consider the calibrating model $\{{U}_{i}  = 
{\beta}_{0}  + {\beta}_{1} {X}_{i}  + {\varepsilon }_{i} ; 1 \leq i \leq
n_{0}\}$ under Gauss--Markov assumptions, such that $\{\QTR{up}{Var}({U}_{i}
) = {\sigma }_{U}^{2} ;1 \leq i \leq n_{0}\}$ and $({\widehat{\beta}}_{0} ,{%
\widehat{\beta}}_{1} )$ are least--squares estimators determining the
empirical calibration line. Under \textit{Gaussian calibration,} the
calibration errors $\{{\varepsilon }_{i} ;1 \leq i \leq n_{0}\}$ comprise 
\textit{iid} ${N}_{1} (0,{\sigma }_{U}^{2} )$ random variables. Subsequent
readings $\{{Z}_{1},\ldots,{Z}_{n} \},$ taken independently of $\{{U}%
_{1},\ldots,{U}_{{n}_{0} } \},$ are then projected as the calibrated
measurements $\{{Y}_{i}  = {\widehat{\beta}}_{0}  + {\widehat{\beta}}_{1} {Z}%
_{i} ;1 \leq i \leq n\}.$ In practice $\{{Z}_{1},\ldots,{Z}_{n} \}$ often
will have been discarded as redundant, or will have been converted directly
without record, so that $\{{Y}_{1},\ldots,{Y}_{n} \}$ remain to be analyzed
and interpreted. If we now suppose that $\boldsymbol{Z} ^{\prime}$ = $[{Z}%
_{1},\ldots,{Z}_{n} ]$ have means $\boldsymbol{\mu} ^{\prime}$ = $[{\mu}%
_{1},\ldots,{\mu}_{n} ]$ and second moments $V(\boldsymbol{Z} )$ = $%
\boldsymbol{\Sigma }  = [{\sigma }_{ij} ],$ independently of $({\widehat{%
\beta}}_{0} ,{\widehat{\beta}}_{1} ),$ then conditional moments of $\mathcal{%
L}({Y}_{1},\ldots,{Y}_{n} \mid{\widehat{\beta}}_{1} )$ are found directly as
follows.

\begin{lemma}
Suppose that $\{{Z}_{1},\ldots,{Z}_{n} \}$ have means $\{{\mu}_{1},\ldots,{%
\mu}_{n} \}$ and second moments $V(\boldsymbol{Z} )$ = $[{\sigma }_{ij} ],$
independently of $({\widehat{\beta}}_{0} ,{\widehat{\beta}}_{1} ),$ and let $%
\{{Y}_{i}  = {\widehat{\beta}}_{0}  + {\widehat{\beta}}_{1} {Z}_{i} ;1 \leq
i \leq n\}.$ Then

(i) $E({Y}_{i} \!\mid\!{\widehat{\beta}}_{1} ) = E({\widehat{\beta}}_{0}
\!\mid\!{\widehat{\beta}}_{1} ) + {\widehat{\beta}}_{1} {\mu}_{i} ;$

(ii) $\QTR{up}{Var}({Y}_{i} \!\mid\!{\widehat{\beta}}_{1} ) = {\widehat{\beta%
}}_{1} ^{2}{\sigma }_{ii}  + \QTR{up}{Var}({\widehat{\beta}}_{0} \!\mid\!{%
\widehat{\beta}}_{1} );$ and

(iii) $\QTR{up}{Cov}({Y}_{i} ,{Y}_{j} \mid{\widehat{\beta}}_{1} )$ = ${%
\widehat{\beta}}_{1} ^{2}{\sigma }_{ij}  + \QTR{up}{Var}({\widehat{\beta}}%
_{0} \!\mid\!{\widehat{\beta}}_{1} ).$
\end{lemma}

\noindent If instead $({\beta}_{0} ,{\beta}_{1} )$ were known, then $\{{Y}%
_{i}  = {\beta}_{0}  + {\beta}_{1} {Z}_{i} ;1 \leq i \leq n\}$ would be
recovered without error, in which case $E({Y}_{i} ) = {\beta}_{0}  + {\beta}%
_{1} {\mu}_{i} $ and $\QTR{up}{Var}({Y}_{i} ) = {\beta}_{1} ^{2}\QTR{up}{Var}%
({Z}_{i} ).$ This ideal case serves as reference against which recovery
subject to calibration errors may be gauged.

Expressions simplify when neither $E({\widehat{\beta}}_{0} \!\mid\!{\widehat{%
\beta}}_{1} )$ nor $\QTR{up}{Var}({\widehat{\beta}}_{0} \!\mid\!{\widehat{%
\beta}}_{1} )$ depends on ${\widehat{\beta}}_{1} ,$ so that $\QTR{up}{Cov}({Y%
}_{i} ,{Y}_{j} \!\mid\!{\widehat{\beta}}_{1} )$ = ${\widehat{\beta}}_{1} ^{2}%
{\sigma }_{ij}  + \QTR{up}{Var}({\widehat{\beta}}_{0} ),$ for example. This
clearly holds under Gaussian calibration, provided that the calibrating
readings $\{{X}_{1},\ldots,{X}_{n_{0}} \}$ have been centered to $\{({X}_{1}
- \overline{X}),\ldots,(X_{n_{0}} - \overline{X})\}.$ This incurs no loss in
generality, as subsequent readings $\{{Z}_{1},\ldots,{Z}_{n} \}$ may be
shifted by $\overline{X}$ units before projecting onto the scale of
measurements. It then follows that ${\widehat{\beta}}_{0}  = \overline{U};$ $%
\QTR{up}{Var}({\widehat{\beta}}_{0} ) = {\sigma }_{0}^{2}  $ = ${\sigma }%
_{U}^{2} /{n}_{0} ;$ and $\QTR{up}{Var}({\widehat{\beta}}_{1} ) = {\sigma }%
_{1}^{2} $ = ${\sigma }_{U}^{2} /{S}_{xx} ,$ where ${S}_{xx} $ = $%
\sum_{i=1}^{n_{0}}({X}_{i}  - \overline{X})^{2};$ so that $({\widehat{\beta}}%
_{0} ,{\widehat{\beta}}_{1} )$ are uncorrelated and thus independent under
Gaussian calibration errors. We henceforth take the initial calibration to
have been centered.

We next consider conditional and unconditional properties of $\{{Y}_{i}  = {%
\widehat{\beta}}_{0}  + {\widehat{\beta}}_{1} {Z}_{i} ;\;1 \leq i \leq n\}$
for the general case that $E(\boldsymbol{Z} )$ = $\boldsymbol{\mu}_{Z} \in%
\mathbb{R}^{n}$ and $V(\boldsymbol{Z} )$ = $\boldsymbol{\Sigma } \in\mathbb{S%
}_{n}^{+},$ to be specialized subsequently. Clearly the conditional means
and dispersion parameters are $E(\boldsymbol{Y} \!\mid\!{\widehat{\beta}}%
_{1} )$ = ${\beta}_{0} \boldsymbol{1}_{n}  + {\widehat{\beta}}_{1} 
\boldsymbol{\mu}_{Z}  = \boldsymbol{\mu}_{Y} ({\widehat{\beta}}_{1} ),$ say,
and $V(\boldsymbol{Y} \!\mid\!{\widehat{\beta}}_{1} ) = {\widehat{\beta}}%
_{1} ^{2}\boldsymbol{\Sigma }  + {\sigma }_{0}^{2} \boldsymbol{1}_{n} 
\boldsymbol{1}_{n} ^{\prime}$ = $\boldsymbol{\Xi} ({\widehat{\beta}}_{1} ).$
Moreover, for the case that $\mathcal{L}(\boldsymbol{Z} )$ = ${N}_{n} (%
\boldsymbol{\mu}_{Z} ,\boldsymbol{\Sigma } )$ in addition to Gaussian
calibration errors, then the conditional distribution of $\boldsymbol{Y} $
is $\mathcal{L}(\boldsymbol{Y} \!\mid\!{\widehat{\beta}}_{1} )$ = ${N}_{n} (%
\boldsymbol{\mu}_{Y} ({\widehat{\beta}}_{1} ),\boldsymbol{\Xi} ({\widehat{%
\beta}}_{1} )).$ Basic unconditional properties follow next.

\begin{theorem}
Consider calibrated measurements $\{{Y}_{i}  = {\widehat{\beta}}_{0}  + {%
\widehat{\beta}}_{1} {Z}_{i} ; 1 \leq i \leq n\}$ projected from readings $\{%
{Z}_{1},\ldots,{Z}_{n} \}$ obtained independently of $({\widehat{\beta}}_{0}
,{\widehat{\beta}}_{1} );$ let $\boldsymbol{Y}  = {\widehat{\beta}}_{0} 
\boldsymbol{1}_{n}  + {\widehat{\beta}}_{1} \boldsymbol{Z} ,$ such that $E(%
\boldsymbol{Z} )$ = $\boldsymbol{\mu}_{Z} \in\mathbb{R}^{n}$ and $V(%
\boldsymbol{Z} )$ = $\boldsymbol{\Sigma } \in\mathbb{S}_{n}^{+};$ and let ${%
\sigma }_{0}^{2}  = \QTR{up}{Var}({\widehat{\beta}}_{0} )$ and ${\sigma }%
_{1}^{2}  = \QTR{up}{Var}({\widehat{\beta}}_{1} ).$ Then unconditional
moments $E(\boldsymbol{Y} ) = \boldsymbol{\mu}_{Y} $ and $V(\boldsymbol{Y} )
= \boldsymbol{\Xi} $ of $\mathcal{L}(\boldsymbol{Y} )$ are given by

(i) $\boldsymbol{\mu}_{Y}  = {\beta}_{0} \boldsymbol{1}_{n}  + {\beta}_{1} 
\boldsymbol{\mu}_{Z} ,$ and

(ii) $\boldsymbol{\Xi}  = {\kappa }_{2} \boldsymbol{\Sigma }  + {\sigma }%
_{0}^{2} \boldsymbol{1}_{n} \boldsymbol{1}_{n} ^{\prime}+ {\sigma }_{1}^{2} 
\boldsymbol{\mu}_{Z} \boldsymbol{\mu}_{Z} ^{\prime},$ with ${\kappa }_{2} $
= ${\mu}_{2} ({\widehat{\beta}}_{1} )$ = ${\sigma }_{1}^{2}  + {\beta}%
_{1}^{2} .$

(iii) Moreover, if $\mathcal{L}(\boldsymbol{Z} )$ = ${N}_{n} (\boldsymbol{\mu%
}_{Z} ,\boldsymbol{\Sigma } )$ in addition to Gaussian calibration errors,
then the unconditional joint density of the elements of $\boldsymbol{Y} $ is
the translation--scale mixture 
\begin{equation}
{f}_{1} (\boldsymbol{y} ;\boldsymbol{\mu}_{Y} ,\boldsymbol{\Xi} ,{G}_{1} ) =
\int_{-\infty}^{\infty} {g}_{n} (\boldsymbol{y} ; \boldsymbol{\mu} (t),%
\boldsymbol{\Xi} (t)) d{G}_{1} (t)
\end{equation}
\noindent as in (2.1), with $\boldsymbol{\mu} (t)$ = ${\beta}_{0} 
\boldsymbol{1}_{n}  + t\boldsymbol{\mu}_{Z} ,$ $\boldsymbol{\Xi} (t) = t^{2}%
\boldsymbol{\Sigma }  + {\sigma }_{0}^{2} \boldsymbol{1}_{n} \boldsymbol{1}%
_{n} ^{\prime},$ and with mixing distribution ${G}_{1} (t) = {N}_{1} ({\beta}%
_{1} , {\sigma }_{1}^{2} ).$
\end{theorem}

\noindent\textbf{Proof.} Conclusion (i) follows directly through
deconditioning. Conclusion (ii) follows using $\QTR{up}{Var}({Y}_{i} )$ = $%
E_{{\widehat{\beta}}_{1} }[\QTR{up}{Var}({Y}_{i} \!\mid\!{\widehat{\beta}}%
_{1} )]$ + $\QTR{up}{Var}_{{\widehat{\beta}}_{1} }(E({Y}_{i} \!\mid\!{%
\widehat{\beta}}_{1} )]$ for variances, and 
\begin{equation}
\QTR{up}{Cov}({Y}_{i} ,{Y}_{j} ) = E_{{\widehat{\beta}}_{1} }[\QTR{up}{Cov}({%
Y}_{i} ,{Y}_{j} \!\mid\!{\widehat{\beta}}_{1} )] + \QTR{up}{Cov}_{{\widehat{%
\beta}}_{1} }[E({Y}_{i} \!\mid\!{\widehat{\beta}}_{1} ),\,E({Y}_{j} \!\mid\!{%
\widehat{\beta}}_{1} )]  \notag
\end{equation}
\noindent for covariances. Conclusion (iii) follows since $\boldsymbol{Y} $
is a linear function of $(\boldsymbol{Z} ,{\widehat{\beta}}_{0} )$ with ${%
\widehat{\beta}}_{1} $ fixed, so that $\mathcal{L}(\boldsymbol{Y} \!\mid\!{%
\widehat{\beta}}_{1} )$ = ${N}_{n} (\boldsymbol{\mu}_{Y} ({\widehat{\beta}}%
_{1} ),\boldsymbol{\Xi} ({\widehat{\beta}}_{1} )),$ as noted, and then
mixing over the distribution of the conditioning variable. $\square$

The foregoing results are basic. We next specialize them as appropriate for
specific experimental settings encountered routinely in practice.

\section{Topics in Inference}

\label{S:topicinf}

Induced dependencies and other model irregularities violate the tenets of
conventional data analysis as noted, specifically, in estimation and
hypothesis testing. We focus on normal--theory procedures, as the
independence typically required by nonparametric competitors is
conspicuously absent. The following sections specialize earlier findings, as
they apply in a single sample, and in one--way comparative experiments.

\subsection{Single Sample}

\label{SS:ssample} Consider $\{{Y}_{i}  = {\widehat{\beta}}_{0}  + {\widehat{%
\beta}}_{1} {Z}_{i} ;1 \leq i \leq n\}$ such that elements of $\boldsymbol{Z}
$ = $[{Z}_{1},\ldots,{Z}_{n} ]^{\prime}$ satisfy $\{E({Z}_{i} ) = {\mu}_{Z}
; 1 \leq i \leq n\}$ and $V(\boldsymbol{Z} )$ = ${\sigma }_{Z}^{2} 
\boldsymbol{I}_{n} .$ If in addition $\{{Z}_{1},\ldots,{Z}_{n} \}$ are 
\textit{iid,} then $\mathcal{L}({Y}_{1},\ldots,{Y}_{n} )$ is exchangeable,
as noted earlier. We are concerned not only with properties of the joint
distribution $\mathcal{L}(\boldsymbol{Y} ),$ but also of $(\overline{Y},{S}%
_{Y}^{2} ,{t}_{0} )$ as the sample mean, the sample variance, and Student's
statistic ${t}_{0} $ = $n^{1/2}(\overline{Y} - {\mu}_{Y}^{0} )/{S}_{Y} ,$ as
well as the ordinary residuals $\{{R}_{i} =({Y}_{i} -\overline{Y});1 \leq i
\leq n\}.$ From Lemma 1 and Theorem 1 we now have the conditional and
unconditional values $\{E({Y}_{i} \!\mid\!{\widehat{\beta}}_{1} ) = {\beta}%
_{0}  + {\widehat{\beta}}_{1} {\mu}_{Z} ;1 \leq i \leq n\};$ $\{E({Y}_{i} )
= {\mu}_{Y}  = {\beta}_{0}  + {\beta}_{1} {\mu}_{Z}  ;1 \leq i \leq n\};$
and $\QTR{up}{V}(\boldsymbol{Y} \!\mid\!{\widehat{\beta}}_{1} )$ = $%
\boldsymbol{\Xi} ({\widehat{\beta}}_{1} )$ = $({\widehat{\beta}}_{1} ^{2}{%
\sigma }_{Z}^{2}  + {\sigma }_{0}^{2} )[(1-\rho)\boldsymbol{I}_{n}  + \rho%
\boldsymbol{1}_{n} \boldsymbol{1}_{n} ^{\prime}],$ with $\rho$ = ${\sigma }%
_{0}^{2} /({\widehat{\beta}}_{1} ^{2}{\sigma }_{Z}^{2} +{\sigma }_{0}^{2} ).$
The unconditional variances are homoscedastic, namely, $\{\QTR{up}{Var}({Y}%
_{i} ) = {\sigma }_{Y}^{2} $ = ${\kappa }_{2} {\sigma }_{Z}^{2}  + {\sigma }%
_{0}^{2}  + {\sigma }_{1}^{2} {\mu}_{Z} ^{2};1 \leq i \leq n\}.$ Essential
findings follow, where it is seen that ${S}_{Y}^{2} $ may grossly
underestimate the actual measurement variance ${\sigma }_{Y}^{2} ,$ and that
structural difficulties becloud both the small--sample and the asymptotic
properties of ${\overline{Y}}_{n} $ = $({Y}_{1} +{Y}_{2} +\ldots+{Y}_{n} )/n.
$

\begin{theorem}
Let $\{{Y}_{i}  = {\widehat{\beta}}_{0}  + {\widehat{\beta}}_{1} {Z}_{i} ; 1
\leq i \leq n\}$ be calibrated measurements from $\{{Z}_{1},\ldots,{Z}_{n}
\},$ such that $E(\boldsymbol{Z} ) = {\mu}_{Z} \boldsymbol{1}_{n} $ and $V(%
\boldsymbol{Z} ) = {\sigma }_{Z}^{2} \boldsymbol{I}_{n} $ independently of $(%
{\widehat{\beta}}_{0} ,{\widehat{\beta}}_{1} );$ and consider the sample
quantities $({\overline{Y}}_{n} ,{S}_{Y}^{2} ),$ together with the ordinary
residuals $\{{R}_{i} =({Y}_{i} -\overline{Y});1 \leq i \leq n\}.$ Then

(i) ${\overline{Y}}_{n} $ is unbiased but inconsistent for estimating ${\mu}%
_{Y} $ = ${\beta}_{0}  + {\beta}_{1} {\mu}_{Z} ;$

(ii) $E({S}_{Y}^{2} ) = {\kappa }_{2} {\sigma }_{Z}^{2} $ = ${\sigma }%
_{Z}^{2} ({\sigma }_{1}^{2}  + {\beta}_{1}^{2} );$ and

(iii) $\{E({R}_{i} ) = 0;1 \leq i \leq n\}.$
\end{theorem}

\noindent\textbf{Proof.} The unbiasedness of ${\overline{Y}}_{n} $ follows
routinely, and its variance from $\QTR{up}{Var}(n^{-1}\boldsymbol{1}_{n}
^{\prime}\boldsymbol{Y} )$ = ${\Xi}_{n} $ with 
\begin{align}
{\Xi}_{n}  &= n^{-2}\boldsymbol{1}_{n} ^{\prime}[{\kappa }_{2} {\sigma }%
_{Z}^{2} \boldsymbol{I}_{n}  + ({\sigma }_{0}^{2}  + {\sigma }_{1}^{2} {\mu}%
_{Z} ^{2})\boldsymbol{1}_{n} \boldsymbol{1}_{n} ^{\prime}]\boldsymbol{1}_{n} 
\\
&= n^{-1}{\kappa }_{2} {\sigma }_{Z}^{2}  + ({\sigma }_{0}^{2}  + {\sigma }%
_{1}^{2} {\mu}_{Z} ^{2}).  \notag
\end{align}
\noindent Since $lim_{n\rightarrow\infty} \QTR{up}{Var}({\overline{Y}}_{n} )$
= $({\sigma }_{0}^{2}  + {\sigma }_{1}^{2} {\mu}_{Z} ^{2}) > 0,$ it follows
that its limit distribution is nondegenerate at ${\mu}_{Y} ,$ so that ${%
\overline{Y}}_{n} $ is consistent neither in probability, nor in mean
square, nor almost surely, as asserted. Conclusion (ii) follows from
evaluating the expected value of the quadratic form $(n-1){S}_{Y}^{2} $ = $%
\boldsymbol{Y} ^{\prime}\boldsymbol{B}_{n} \boldsymbol{Y} $ as $E[(n-1){S}%
_{Y}^{2} ]$ = $\QTR{up}{tr}\,\boldsymbol{B}_{n}  V(\boldsymbol{Y} )$ + $%
\boldsymbol{\mu}_{Y} ^{\prime}\boldsymbol{B}_{n} \boldsymbol{\mu}_{Y} .$
Details are 
\begin{align}
E[(n-1){S}_{Y}^{2} ] &= \QTR{up}{tr}\,\boldsymbol{B}_{n} [{\kappa }_{2} {%
\sigma }_{Z}^{2} \boldsymbol{I}_{n}  + ({\sigma }_{0}^{2}  + {\sigma }%
_{1}^{2} {\mu}_{Z} ^{2})\boldsymbol{1}_{n} \boldsymbol{1}_{n} ^{\prime}]%
\boldsymbol{B}_{n}  + \boldsymbol{\mu}_{Y} ^{\prime}\boldsymbol{B}_{n} 
\boldsymbol{\mu}_{Y}   \notag \\
&= (n-1){\kappa }_{2} {\sigma }_{Z}^{2}   \notag
\end{align}
\noindent where $\boldsymbol{\mu}_{Y} ^{\prime}\boldsymbol{B}_{n} 
\boldsymbol{\mu}_{Y} $ = $({\beta}_{0}  + {\beta}_{1} {\mu}_{Z} )^{2}%
\boldsymbol{1}_{n} ^{\prime}\boldsymbol{B}_{n} \boldsymbol{1}_{n} $ = 0,
since $\boldsymbol{B}_{n} $ is idempotent and $\boldsymbol{B}_{n} 
\boldsymbol{1}_{n}  = \boldsymbol{0} .$ Conclusion (iii) follows directly,
to complete our proof. $\square$

The following consequences may be noted.

\begin{itemize}
\item Conclusion (i) appears to dash the usual expectation that lengths of $%
(1-\alpha )$ confidence intervals for ${\mu}_{Y} $ will decrease at the rate 
$O(n^{-1/2}).$

\item The sample variance ${S}_{Y}^{2} $ underestimates the actual variance $%
{\sigma }_{Y}^{2} .$ The bias is $E({S}_{Y}^{2} )$ - ${\sigma }_{Y}^{2} $ = $%
- ({\sigma }_{0}^{2}  + {\sigma }_{1}^{2} {\mu}_{Z} ^{2}).$

\item This bias increases with decreasing precision in estimating the
calibration line, and with increasing $\mid\!{\mu}_{Z} \!\mid$ and thus $%
\mid\!{\mu}_{Y} \!\mid.$

\item On the other hand, the expectation $E({S}_{Y}^{2} )$ = ${\kappa }_{2} {%
\sigma }_{Z}^{2} ,$ with ${\kappa }_{2} $ = ${\mu}_{2} ({\widehat{\beta}}%
_{1} ),$ may be compared with the ideal variance, $\QTR{up}{Var}({Y}_{i} ) = 
{\beta}_{1} ^{2}{\sigma }_{Z}^{2} ,$ attained under linear calibration with
known parameters $({\beta}_{0} ,{\beta}_{1} ).$
\end{itemize}

We have seen how unconditional moments of calibrated measurements depend on
those of the conditioning variable ${\widehat{\beta}}_{1} .$ It remains to
examine effects of the fitted calibration line on unconditional
distributions, to include those of various sample statistics. Recall from
Theorem 2 and its proof that $E(\overline{Y})$ = ${\mu}_{Y} $ = ${\beta}_{0}
+ {\beta}_{1} {\mu}_{Z} $ and $\QTR{up}{Var}(\overline{Y})$ = $n^{-1}{\kappa 
}_{2} {\sigma }_{Z}^{2}  + ({\sigma }_{0}^{2}  + {\sigma }_{1}^{2} {\mu}_{Z}
^{2}).$ To invoke expression (2.1) and its special case at $\boldsymbol{%
\theta}  = \boldsymbol{0} ,$ under Gaussian assumptions we have ${G}_{1} ({%
\widehat{\beta}}_{1} )$ = ${N}_{1} ({\beta}_{1} ,{\sigma }_{1}^{2} ),$
together with ${G}_{2} ({\widehat{\beta}}_{1} ^{2};\lambda ),$ such that $%
\mathcal{L}({\widehat{\beta}}_{1} ^{2}/{\sigma }_{1}^{2} )$ = ${\chi}^{2}
(1,\lambda )$ with $\lambda  = {\beta}_{1} ^{2}/{\sigma }_{1}^{2} .$ Basic
unconditional distributions follow next as mixtures .

\begin{theorem}
Let $\{{Y}_{i}  = {\widehat{\beta}}_{0}  + {\widehat{\beta}}_{1} {Z}_{i} ; 1
\leq i \leq n\}$ be calibrated measurements; suppose that $\{{Z}_{1},\ldots,{%
Z}_{n} \}$ are \textit{iid} ${N}_{1} ({\mu}_{Z} ,{\sigma }_{Z}^{2} )$
independently of $({\widehat{\beta}}_{0} ,{\widehat{\beta}}_{1} )$ under
Gaussian calibration; and consider the sample quantities $({\overline{Y}}%
_{n} ,{S}_{Y}^{2} ,{t}_{0}^{2} ),$ together with the ordinary residuals $\{{R%
}_{i} =({Y}_{i} -\overline{Y});1 \leq i \leq n\},$ where ${t}_{0}^{2}  = n(%
\overline{Y} - {\mu}_{Y}^{0} )^{2}/{S}_{Y}^{2} $ for testing $H_{0}:{\mu}%
_{Y}  = {\mu}_{Y}^{0} $ against $H_{1}:{\mu}_{Y}  \neq {\mu}_{Y}^{0} .$ Then
unconditional properties are as follows.

(i) The unconditional density of $\mathcal{L}(\overline{Y})$ is the
translation--scale mixture 
\begin{equation}
{f}_{1} (u;{\mu}_{Y} ,{\Xi}_{n} ,{G}_{1} ) = \int_{-\infty}^{\infty} {g}_{1}
(u; \mu(t),{\Xi}_{n} (t))d{G}_{1} (t)
\end{equation}
\noindent with mixing distribution ${G}_{1} ({\widehat{\beta}}_{1} ) = {N}%
_{1} ({\beta}_{1} , {\sigma }_{1}^{2} )$ as in (2.1) for the case $n=1,$
where $\mu(t) = {\beta}_{0}  + t{\mu}_{Z} $ and ${\Xi}_{n} (t) = (t^{2}{%
\sigma }_{Z}^{2} /n + {\sigma }_{0}^{2} ),$ and ${\Xi}_{n} $ is defined in
(4.1).

(ii) The joint density of residuals $\boldsymbol{R} $ = $[{R}_{1},\ldots,{R}%
_{n} ]^{\prime}$ is given by 
\begin{equation}
{f}_{n} (\boldsymbol{r} ;\boldsymbol{0} ,{\sigma }_{Z}^{2} \boldsymbol{B}%
_{n} ,{G}_{2} ) = \int_{0}^{\infty}{g}_{n} (\boldsymbol{r} ;\boldsymbol{0} ,t%
{\sigma }_{Z}^{2} \boldsymbol{B}_{n} )d{G}_{2} (t)
\end{equation}
\noindent with mixing distribution ${G}_{2} ({\widehat{\beta}}_{1}
^{2};\lambda )$ = ${\chi}^{2} (1,\lambda )$ and $\lambda  = {\beta}_{1} ^{2}/%
{\sigma }_{1}^{2} .$

(iii) The joint distribution $\mathcal{L}({R}_{1},\ldots,{R}_{n} )$
increases in peakedness about $\boldsymbol{0} \in\mathbb{R}^{n}$ with
decreasing $\lambda .$

(iv) $\mathcal{L}(\nu{S}_{Y}^{2} /{\sigma }_{Z}^{2} )$ has the density ${f}%
_{0} (u;\nu/2,2,{G}_{2} )$ as in (2.2), with $\nu = n-1,$ mixing
distribution ${G}_{2} ({\widehat{\beta}}_{1} ^{2};\lambda )$ = ${\chi}^{2}
(1,\lambda )$ and $\lambda  = {\beta}_{1}^{2} /{\sigma }_{1}^{2} ;$ and $E({S%
}_{Y}^{2} )$ = ${\sigma }_{Z}^{2} ({\sigma }_{1}^{2}  + {\beta}_{1}^{2} ).$

(v) The distribution $\mathcal{L}(\nu{S}_{Y}^{2} /{\sigma }_{Z}^{2} )$
increases stochastically with $\lambda  = {\beta}_{1} ^{2}/{\sigma }_{1}^{2}
.$

(vi) The unconditional density of ${t}_{0}^{2}  = n(\overline{Y} - {\mu}%
_{Y}^{0} )^{2}/{S}_{Y}^{2} $ is given as the mixture 
\begin{equation}
g(u;\nu,\delta ,{G}_{2} ) = \int_{0}^{\infty} {g}_{T^{2}} (u;\nu,\delta /t)d{%
G}_{2} (t)  \notag
\end{equation}
\noindent with mixing distribution ${G}_{2} ({\widehat{\beta}}_{1}
^{2};\lambda ),$ where $\nu = n-1,$ $\delta  = ({\mu}_{Y}  - {\mu}_{Y}^{0}
)^{2}/{\sigma }_{Z}^{2} ,$ and $\lambda $ = ${\beta}_{1} ^{2}/{\sigma }%
_{1}^{2} .$

(vii) The unconditional \textit{cdf} of $\mathcal{L}({t}_{0}^{2} )$
increases stochastically with increasing $\delta  = ({\mu}_{Y}  - {\mu}%
_{Y}^{0} )^{2}/{\sigma }_{Z}^{2} $ for fixed $\lambda $ = ${\beta}_{1} ^{2}/{%
\sigma }_{1}^{2} ;$ and for fixed $\delta ,$ it decreases stochastically
with increasing $\lambda .$
\end{theorem}

\noindent\textbf{Proof.} The conditional distribution of note is $\mathcal{L}%
(\overline{Y}\!\mid\!{\widehat{\beta}}_{1} )$ = ${N}_{1} (\mu({\widehat{\beta%
}}_{1} ),{\Xi}_{n} ({\widehat{\beta}}_{1} )),$ where $\mu({\widehat{\beta}}%
_{1} ) = {\beta}_{0}  + {\widehat{\beta}}_{1} {\mu}_{Z} $ and ${\Xi}_{n} ({%
\widehat{\beta}}_{1} )$ = $({\widehat{\beta}}_{1} ^{2}{\sigma }_{Z}^{2} /n + 
{\sigma }_{0}^{2} ).$ Its unconditional density thus is ${f}_{1} (u;{\mu}%
_{Y} ,{\Xi}_{n} ,{G}_{1} )$ from (2.1), to give conclusion (i) with mixing
distribution as asserted. To continue, observe that $\boldsymbol{R} $ = $%
\boldsymbol{B}_{n} \boldsymbol{Y} $ and $\boldsymbol{B}_{n} \boldsymbol{\Xi}
({\widehat{\beta}}_{1} )\boldsymbol{B}_{n} $ = $\boldsymbol{B}_{n} ({%
\widehat{\beta}}_{1} ^{2}{\sigma }_{Z}^{2} \boldsymbol{I}_{n}  + {\sigma }%
_{0}^{2} \boldsymbol{1}_{n} \boldsymbol{1}_{n} ^{\prime})\boldsymbol{B}_{n} 
= {\widehat{\beta}}_{1} ^{2}{\sigma }_{Z}^{2} \boldsymbol{B}_{n} $ since $%
\boldsymbol{B}_{n} $ is idempotent and $\boldsymbol{B}_{n} \boldsymbol{1}%
_{n}  = \boldsymbol{0} .$ We infer conditionally that $\mathcal{L}(%
\boldsymbol{R} \!\mid\!{\widehat{\beta}}_{1} )$ = ${N}_{n} (\boldsymbol{0} ,{%
\widehat{\beta}}_{1} ^{2}{\sigma }_{Z}^{2} \boldsymbol{B}_{n} ),$ since $%
\boldsymbol{B}_{n} \boldsymbol{\mu} ({\widehat{\beta}}_{1} )$ = $({\beta}%
_{0}  + {\widehat{\beta}}_{1} {\mu}_{Z} )\boldsymbol{B}_{n} \boldsymbol{1}%
_{n}  = \boldsymbol{0} .$ The unconditional distribution is the scale
mixture as in conclusion (ii), with ${G}_{2} ({\widehat{\beta}}_{1}
^{2},\lambda )$ = ${\chi}^{2} (1,\lambda )$ as the mixing distribution over $%
[0,\infty).$ Conclusion (iii) follows from Lemma A2 of the Appendix, as the
mixing distribution $\mathcal{L}({\widehat{\beta}}_{1} ^{2}/{\sigma }%
_{1}^{2} )$ = ${\chi}^{2} (1,\lambda )$ increases stochastically with $%
\lambda  = {\beta}_{1} ^{2}/{\sigma }_{1}^{2} .$ To proceed, observe that $%
(n-1){S}_{Y}^{2} $ = $\boldsymbol{R} ^{\prime}\boldsymbol{R} ,$ so that $%
\mathcal{L}(\boldsymbol{R} ^{\prime}\boldsymbol{R} /{\widehat{\beta}}_{1}
^{2}{\sigma }_{Z}^{2} \!\mid\!{\widehat{\beta}}_{1} ^{2})$ = ${\chi}^{2}
(n-1,0).$ It follows that $\mathcal{L}[(n-1){S}_{Y}^{2} /{\sigma }_{Z}^{2}
\!\mid\!{\widehat{\beta}}_{1} ^{2}]$ is a central chi--squared variate
scaled by ${\widehat{\beta}}_{1} ^{2}.$ On identifying $(n-1){S}_{Y}^{2} /{%
\sigma }_{Z}^{2} $ with $U$ and ${\widehat{\beta}}_{1} ^{2}$ with $w$ in
developments leading to (2.2), we thus establish conclusion (iv) on
specializing from gamma to chi--squared distributions. Conclusion (v)
follows directly from (iii) since the set ${A}_{t} $ = $\{\boldsymbol{R} \in%
\mathbb{R}^{n}\,:\,\boldsymbol{R} ^{\prime}\boldsymbol{R}  \leq t\}$ is
convex and symmetric in $\mathbb{R}^{n},$ whereas $P({A}_{t} )$ = $P[(n-1){S}%
_{Y}^{2} \leq t{\sigma }_{Z}^{2} ].$

To see conclusion (vi), recall the affine--invariance of ${t}^{2}  = n(%
\overline{Y} - {\mu}_{Y}^{0} )^{2}/{S}_{Y}^{2} $ for testing $H_{0}:\,{\mu}%
_{Y}  = {\mu}_{Y}^{0} $ under $\{{Y}_{i} \to a + b{Z}_{i} ;\,1 \leq i \leq
n\},$ namely, ${t}^{2}  = n[(\overline{Y} - {\mu}_{Y} ) + ({\mu}_{Y}  - {\mu}%
_{Y}^{0} )]^{2}/{S}_{Y}^{2} $ = $n\{[(a + b\overline{Z}) - (a + b{\mu}_{Z}
)] + b({\mu}_{Y}  - {\mu}_{Y}^{0} )/b\}^{2}/b^{2}{S}_{Z}^{2} $ = $n[(%
\overline{Z} - {\mu}_{Z} ) + ({\mu}_{Y}  - {\mu}_{Y}^{0} )/b]^{2}/{S}%
_{Z}^{2} .$ Clearly $\mathcal{L}({t}^{2} )$ = ${t}^{2} (\nu,\delta ),$ with $%
\nu = n-1$ and $\delta  = ({\mu}_{Y}  - {\mu}_{Y}^{0} )^{2}/b^{2}{\sigma }%
_{Z}^{2} ,$ independently of $a.$ We next apply these facts conditionally,
given $({\widehat{\beta}}_{0} ,{\widehat{\beta}}_{1} ),$ on replacing $(a,b)$
with $({\widehat{\beta}}_{0} ,{\widehat{\beta}}_{1} ),$ to infer that $%
\mathcal{L}({t}_{0}^{2} \!\mid\!{\widehat{\beta}}_{0} ,{\widehat{\beta}}_{1}
)$ = $\mathcal{L}({t}_{0}^{2} \!\mid\!{\widehat{\beta}}_{1} )$ = ${t}^{2}
(\nu,\delta ({\widehat{\beta}}_{1} )),$ where $\delta ({\widehat{\beta}}_{1}
) = ({\mu}_{Y}  - {\mu}_{Y}^{0} )^{2}/{\widehat{\beta}}_{1} ^{2}{\sigma }%
_{Z}^{2} ,$ independently of ${\widehat{\beta}}_{0} .$ It follows that the
unconditional distribution of ${t}_{0}^{2} $ has the mixture density $%
g(u;\nu,\delta ,{G}_{2} )$ = $\int_{0}^{\infty} {g}_{T^{2}} (u;\nu,\delta
/u)d{G}_{2} (u)$ with mixing distribution ${G}_{2} ({\widehat{\beta}}_{1}
^{2};\lambda )$ as asserted, where $\delta  = ({\mu}_{Y}  - {\mu}_{Y}^{0}
)^{2}/{\sigma }_{Z}^{2} $ and $\lambda  = {\beta}_{1} ^{2}/{\sigma }_{1}^{2}
.$ To continue, the test for $H_{0}:{\mu}_{Y}  = {\mu}_{Y}^{0} $ against $%
H_{1}:{\mu}_{Y}  \neq {\mu}_{Y}^{0} $ rejects at level $\alpha $ for ${t}%
_{0}^{2}  > {c}_{\alpha }^{2} ;$ moreover, the conditional \textit{cdf} $%
\mathcal{L}({t}_{0}^{2} \!\mid\!{\widehat{\beta}}_{1} )$ increases
stochastically with $\delta ({\widehat{\beta}}_{1} ),$ point--wise for each
fixed ${\widehat{\beta}}_{1} ,$ from standard properties of noncentral ${t}%
^{2} $--distributions. It follows that the unconditional \textit{cdf}
increases stochastically with increasing $\delta $ under mixing. That the
unconditional \textit{cdf} decreases stochastically with increasing $\lambda
,$ with $\delta $ held fixed, follows unconditionally through mixing as in
the proof for conclusion (v), to complete our proof. $\square$

Note from conclusion (iv) that $E({S}_{Y}^{2} ) = {\sigma }_{Z}^{2} ({\sigma 
}_{1}^{2}  + {\beta}_{1}^{2} ).$ This may be compared with the ideal case $%
\QTR{up}{Var}({Y}_{i} ) = {\beta}_{1}^{2} {\sigma }_{Z}^{2} $ where $({\beta}%
_{0} ,{\beta}_{1} )$ are known. It is instructive to reexamine unconditional
properties of this section if we retain the homogeneity of variances of $\{{Z%
}_{1},\ldots,{Z}_{n} \},$ but assume instead that means are not, \textit{%
i.e.,} that $E(\boldsymbol{Z} )$ = $\boldsymbol{\mu}  = [{\mu}_{1},\ldots,{%
\mu}_{n} ]^{\prime}.$ For then we see that $\{\QTR{up}{Var}({Y}_{i} )$ = ${%
\kappa }_{2} {\sigma }_{Z}^{2}  + {\sigma }_{0}^{2}  + {\sigma }_{1}^{2} {\mu%
}_{i} ^{2};1 \leq i \leq n\}.$ We have the curious finding that calibrated
measurements in a single sample will have heterogeneous variances under
heterogeneous means, despite the homogeneity of variances of $\{{Z}%
_{1},\ldots,{Z}_{n} \}.$

\subsection{One--Way Experiments}

\label{SS:oneway} Clearly $\{{Y}_{1},\ldots,{Y}_{n} \}$ and $\{{Z}%
_{1},\ldots,{Z}_{n} \}$ have the same experimental structure, here a
one--way experiment comprising $k$ samples of sizes $\{{n}_{1},\ldots,{n}%
_{k} \},$ with ${n}_{1} + \ldots + {n}_{k}  = n.$ In keeping with
conventional notation, partition $\boldsymbol{Z} ^{\prime}= [{Z}_{1},\ldots,{%
Z}_{n} ]$ as $[\boldsymbol{Z}_{1} ^{\prime},\ldots,\boldsymbol{Z}_{k}
^{\prime}],$ such that $\{\boldsymbol{Z}_{i} ^{\prime}= [{Z}_{i1},\ldots,{Z}%
_{i{n}_{i} } ];1 \leq i \leq k\};$ and similarly for $\boldsymbol{Y}
^{\prime}$ = $[\boldsymbol{Y}_{1} ^{\prime},\ldots,\boldsymbol{Y}_{k}
^{\prime}],$ with $\{\boldsymbol{Y}_{i} ^{\prime}= [{Y}_{i1},\ldots,{Y}_{i{n}%
_{i} } ];1 \leq i \leq k\};$ and suppose that $\{E({Z}_{ij} ) = {\mu}_{i} ;1
\leq j \leq {n}_{i} \}$ and $\{\QTR{up}{Var}({Z}_{ij} ) = {\omega }_{i}^{2}
;1 \leq j \leq {n}_{i} \}.$ Accordingly, take $\boldsymbol{\mu}_{Z} $ = $[{%
\mu}_{1} \boldsymbol{1}_{{n}_{1} } ^{\prime},\ldots,{\mu}_{k} \boldsymbol{1}%
_{{n}_{k} } ^{\prime}]^{\prime}$ and $V(\boldsymbol{Z} )$ = $\QTR{up}{Diag}({%
\omega }_{1}^{2} \boldsymbol{I}_{{n}_{1} } ,\ldots,{\omega }_{k}^{2} 
\boldsymbol{I}_{{n}_{k} } )$ = $\boldsymbol{D} ({\omega }_{1}^{2},\ldots,{%
\omega }_{k}^{2} ),$ say. Specializing from Section 3, we have the
conditional moments $E(\boldsymbol{Y} \!\mid\!{\widehat{\beta}}_{1} )$ = $%
\boldsymbol{\mu}_{Y} ({\widehat{\beta}}_{1} )$ = ${\beta}_{0} \boldsymbol{1}%
_{n}  + {\widehat{\beta}}_{1}  [{\mu}_{1} \boldsymbol{1}_{{n}_{1} }
^{\prime},\ldots,{\mu}_{k} \boldsymbol{1}_{{n}_{k} } ^{\prime}]^{\prime}$
and $V(\boldsymbol{Y} \!\mid\!{\widehat{\beta}}_{1} )$ = $\boldsymbol{\Xi} ({%
\widehat{\beta}}_{1} )$ = ${\widehat{\beta}}_{1} ^{2}\QTR{up}{Diag}({\omega }%
_{1}^{2} \boldsymbol{I}_{{n}_{1} } ,\ldots,{\omega }_{k}^{2} \boldsymbol{I}_{%
{n}_{k} } ) + {\sigma }_{0}^{2} \boldsymbol{1}_{n} \boldsymbol{1}_{n}
^{\prime},$ together with $\mathcal{L}(\boldsymbol{Y} \!\mid\!{\widehat{\beta%
}}_{1} )$ = ${N}_{n} (\boldsymbol{\mu}_{Y} ({\widehat{\beta}}_{1} ),%
\boldsymbol{\Xi} ({\widehat{\beta}}_{1} ))$ under Gaussian errors. Moreover,
unconditional moments are $E(\boldsymbol{Y} )$ = $\boldsymbol{\mu}_{Y} $ = ${%
\beta}_{0} \boldsymbol{1}_{n}  + {\beta}_{1}  [{\mu}_{1} \boldsymbol{1}_{{n}%
_{1} } ^{\prime},\ldots,{\mu}_{k} \boldsymbol{1}_{{n}_{k} }
^{\prime}]^{\prime}$ and $V(\boldsymbol{Y} )$ = $\boldsymbol{\Xi} $ = ${%
\kappa }_{2} \QTR{up}{Diag}({\omega }_{1}^{2} \boldsymbol{I}_{{n}_{1} }
,\ldots,{\omega }_{k}^{2} \boldsymbol{I}_{{n}_{k} } )$ + ${\sigma }_{0}^{2} 
\boldsymbol{1}_{n} \boldsymbol{1}_{n} ^{\prime}$ + ${\sigma }_{1}^{2} 
\boldsymbol{M} ,$ where $\boldsymbol{M}  = [\boldsymbol{M}_{ij} ]$ = $%
\boldsymbol{\mu}_{Z} \boldsymbol{\mu}_{Z} ^{\prime}$ has the block structure 
$\boldsymbol{M}_{ij}  = {\mu}_{i} {\mu}_{j} \boldsymbol{1}_{{n}_{i} } 
\boldsymbol{1}_{{n}_{j} } ^{\prime}.$ In particular, for typical calibrated
measurements in sample $i$ of the $k$ samples, the conditional and
unconditional means are $\{E({Y}_{ij} \!\mid\!{\widehat{\beta}}_{1} ) = {%
\beta}_{0}  + {\widehat{\beta}}_{1} {\mu}_{i} ;1 \leq j \leq {n}_{i} \}$ and 
$\{E({Y}_{ij} ) = {\beta}_{0}  + {\beta}_{1} {\mu}_{i} ;1 \leq j \leq {n}%
_{i} \},$ whereas the corresponding variances are $\{\QTR{up}{Var}({Y}_{ij}
\!\mid\!{\widehat{\beta}}_{1} )$ = ${\widehat{\beta}}_{1} ^{2}{\omega }%
_{i}^{2}  + {\sigma }_{0}^{2} ;1 \leq j \leq {n}_{i} \}$ and $\{\QTR{up}{Var}%
({Y}_{ij} ) = {\kappa }_{2} {\omega }_{i}^{2}  + {\sigma }_{0}^{2}  + {%
\sigma }_{1}^{2} {\mu}_{i} ^{2};1 \leq j \leq {n}_{i} \}.$ We are concerned
with the dual issues of the homogeneity of means, and of the homogeneity of
variances, across the $k$ samples. Clearly the induced irregularities are
artifacts of the calibration process, rather than consequences of the
experimental structure itself. It remains to determine unintended effects of
calibration on conventional comparisons among the means and the variances.

In contrast to conventional one--way experiments, where homoscedasticity can
be checked regardless of heterogeneity among the $k$ population means, under
calibration we see that homogeneity of the unconditional variances is
possible only in unusual circumstances. Specifically, homoscedasticity holds
unconditionally if and only if, for every pair $({\omega }_{i}^{2} ,{\mu}%
_{i}^{2} )$ and $({\omega }_{j}^{2} ,{\mu}_{j}^{2} ),$ that $({\omega }%
_{i}^{2}  - {\omega }_{j}^{2} )$ = $c({\mu}_{j}^{2}  - {\mu}_{i}^{2} )$ with 
$c$ = ${\sigma }_{1}^{2} /{\kappa }_{2} .$

To continue, consider transformations ${T}_{1} (\boldsymbol{Y} ),$ ${T}_{2} (%
\boldsymbol{Y} ),$ and ${T}_{3} (\boldsymbol{Y} )$ such that ${T}_{1} (%
\boldsymbol{Y} )$ = $\bar{\boldsymbol{Y} }$ = $[{\overline{Y}}_{1},\ldots,{%
\overline{Y}}_{k} ]^{\prime}$ comprise the $k$ sample means; ${T}_{2} (%
\boldsymbol{Y} )$ = $\boldsymbol{R} ^{\prime}$ = $[\boldsymbol{R}_{1}
^{\prime},\ldots,\boldsymbol{R}_{k} ^{\prime}]^{\prime}$ consists of the
ordinary within--sample residuals, with $\boldsymbol{R}_{i} $ = $\boldsymbol{%
B}_{n_{i}} \boldsymbol{Y}_{i} $ and $\boldsymbol{B}_{n_{i}}  = (\boldsymbol{I%
}_{{n}_{i} }  - {n}_{i} ^{-1}\boldsymbol{1}_{{n}_{i} } \boldsymbol{1}_{{n}%
_{i} } ^{\prime});$ and ${T}_{3} (\boldsymbol{Y} )$ = $[{\nu}_{1} {S}%
_{1}^{2} ,\ldots,{\nu}_{k} {S}_{k}^{2} ]^{\prime}$ are the residual sums of
squares, \textit{i.e.,} $\{{\nu}_{i} {S}_{i}^{2}  = \boldsymbol{R}_{i}
^{\prime}\boldsymbol{R}_{i}  = \boldsymbol{Y}_{i} ^{\prime}\boldsymbol{B}%
_{n_{i}} \boldsymbol{Y}_{i} ;1 \leq i \leq k\},$ with ${\nu}_{i}  = {n}_{i}
-1.$ Basic properties may be summarized as follows.

\begin{theorem}
Consider calibrated measurements $\boldsymbol{Y} ^{\prime}= [\boldsymbol{Y}%
_{1} ^{\prime},\ldots,\boldsymbol{Y}_{k} ^{\prime}]$ corresponding to $%
\boldsymbol{Z} ^{\prime}= [\boldsymbol{Z}_{1} ^{\prime},\ldots,\boldsymbol{Z}%
_{k} ^{\prime}]$ such that $E(\boldsymbol{Z} )$ = $[{\mu}_{1} \boldsymbol{1}%
_{{n}_{1} } ^{\prime},\ldots,{\mu}_{k} \boldsymbol{1}_{{n}_{k} }
^{\prime}]^{\prime}$ and $V(\boldsymbol{Z} )$ = $\QTR{up}{Diag}({\omega }%
_{1}^{2} \boldsymbol{I}_{{n}_{1} } ,\ldots,{\omega }_{k}^{2} \boldsymbol{I}_{%
{n}_{k} } );$ and let ${T}_{1} (\boldsymbol{Y} )$ = $\bar{\boldsymbol{Y} }$
= $[{\overline{Y}}_{1},\ldots,{\overline{Y}}_{k} ]^{\prime},$ ${T}_{2} (%
\boldsymbol{Y} )$ = $[\boldsymbol{R}_{1} ^{\prime},\ldots,\boldsymbol{R}_{k}
^{\prime}]^{\prime},$ and ${T}_{3} (\boldsymbol{Y} )$ = $[{\nu}_{1} {S}%
_{1}^{2} ,\ldots,{\nu}_{k} {S}_{k}^{2} ]^{\prime},$ with $\{{\nu}_{i}  = {n}%
_{i} -1;1 \leq i \leq k\}.$ Moreover, a Gaussian model asserts that $\{({Z}%
_{ij}  - {\mu}_{i} )/{\omega }_{i}^{2} ; 1 \leq j \leq {n}_{i} ,1 \leq i
\leq k\}$ are \textit{iid} ${N}_{1} (0,1)$ random variates independently of $%
({\widehat{\beta}}_{0} ,{\widehat{\beta}}_{1} )$ under Gaussian calibration.

(i) Conditional and unconditional moments of ${T}_{1} (\boldsymbol{Y} ) = 
\bar{\boldsymbol{Y} }$ are given by $E(\bar{\boldsymbol{Y} }\!\mid\!{%
\widehat{\beta}}_{1} )$ = $\boldsymbol{\theta} ({\widehat{\beta}}_{1} )$ = ${%
\beta}_{0} \boldsymbol{1}_{k}  + {\widehat{\beta}}_{1} \boldsymbol{\mu} ,$ $%
E(\bar{\boldsymbol{Y} })$ = $\boldsymbol{\theta} $ = ${\beta}_{0} 
\boldsymbol{1}_{k}  + {\beta}_{1} \boldsymbol{\mu} ,$ $V(\bar{\boldsymbol{Y} 
}\!\mid\!{\widehat{\beta}}_{1} )$ = $\boldsymbol{\Xi}_{1} ({\widehat{\beta}}%
_{1} )$ = ${\widehat{\beta}}_{1} ^{2}\QTR{up}{Diag}({\omega }_{1}^{2} /{n}%
_{1} ,\ldots,{\omega }_{k}^{2} /{n}_{k} ) +$ ${\sigma }_{0}^{2} \boldsymbol{1%
}_{k} \boldsymbol{1}_{k} ^{\prime},$ and $V(\bar{\boldsymbol{Y} })$ = $%
\boldsymbol{\Xi}_{1} $ = ${\kappa }_{2} \QTR{up}{Diag}({\omega }_{1}^{2} /{n}%
_{1} ,\ldots,{\omega }_{k}^{2} /{n}_{k} ) + {\sigma }_{0}^{2} \boldsymbol{1}%
_{k} \boldsymbol{1}_{k} ^{\prime}+ {\sigma }_{1}^{2} \boldsymbol{\mu} 
\boldsymbol{\mu} ^{\prime},$ where $\boldsymbol{\mu} ^{\prime}$ = $[{\mu}%
_{1},\ldots,{\mu}_{k} ].$

(ii) Under Gaussian models the unconditional density of $\mathcal{L}(\bar{%
\boldsymbol{Y} })$ is the translation--scale mixture 
\begin{equation}
{f}_{k} (\boldsymbol{u} ;\boldsymbol{\theta} ,\boldsymbol{\Xi}_{1} ,{G}_{1}
) = \int_{-\infty}^{\infty} {g}_{k} (\boldsymbol{u} ; \boldsymbol{\theta}
(t),\boldsymbol{\Xi}_{1} (t))d{G}_{1} (t)
\end{equation}
\noindent with mixing distribution ${G}_{1} ({\widehat{\beta}}_{1} ) = {N}%
_{1} ({\beta}_{1} , {\sigma }_{1}^{2} )$ as in (2.1), where $\boldsymbol{%
\theta} (t) = {\beta}_{0} \boldsymbol{1}_{k}  + t\boldsymbol{\mu} $ and $%
\boldsymbol{\Xi}_{1} (t)$ = $t^{2}\QTR{up}{Diag}({\omega }_{1}^{2} /{n}_{1}
,\ldots,{\omega }_{k}^{2} /{n}_{k} ) + {\sigma }_{0}^{2} \boldsymbol{1}_{k} 
\boldsymbol{1}_{k} ^{\prime}.$

(iii) Conditional and unconditional moments of the residuals are $E(%
\boldsymbol{R} \!\mid\!{\widehat{\beta}}_{1} )$ = $E(\boldsymbol{R} ) = 
\boldsymbol{0} ,$ $V(\boldsymbol{R} \!\mid\!{\widehat{\beta}}_{1} )$ = $%
\boldsymbol{\Xi}_{2} ({\widehat{\beta}}_{1} )$ = ${\widehat{\beta}}_{1} ^{2}%
\QTR{up}{Diag}({\omega }_{1}^{2} \boldsymbol{B}_{n_{1}} ,\ldots,{\omega }%
_{k}^{2} \boldsymbol{B}_{n_{k}} ),$ and 
\begin{equation}
V(\boldsymbol{R} ) = \boldsymbol{\Xi}_{2}  = {\kappa }_{2} \QTR{up}{Diag}({%
\omega }_{1}^{2} \boldsymbol{B}_{n_{1}} ,\ldots,{\omega }_{k}^{2} 
\boldsymbol{B}_{n_{k}} ).  \notag
\end{equation}

(iv) Under Gaussian models the joint density of residuals $\boldsymbol{R} $
= $[\boldsymbol{R}_{1} ^{\prime},\ldots,\boldsymbol{R}_{k}
^{\prime}]^{\prime}$ is given by ${f}_{2} (\boldsymbol{r} ;\boldsymbol{0} ,%
\boldsymbol{\Xi}_{2} ,{G}_{2} )$ as in (4.3), with mixing distribution ${G}%
_{2} ({\widehat{\beta}}_{1} ^{2};\lambda )$ and $\lambda  = {\beta}_{1} ^{2}/%
{\sigma }_{1}^{2} .$

(v) Under Gaussian models the joint density of elements of $[{\nu}_{1} {S}%
_{1}^{2} /{\omega }_{1}^{2} ,\ldots,{\nu}_{k} {S}_{k}^{2} /{\omega }_{k}^{2}
]^{\prime}$ is given by 
\begin{equation}
f(\boldsymbol{u} ;{\nu}_{1},\ldots,{\nu}_{k} ) = \int_{0}^{\infty}
\prod_{i=1}^{k} {g}_{0} ({u}_{i} ;{\nu}_{i} /2,2w)d{G}_{2} (w)  \notag
\end{equation}
\noindent with ${\nu}_{i}  = {n}_{i} -1$ and ${g}_{0} (u;\alpha ,w\beta )$ = 
$(w\beta )^{-\alpha } u^{\alpha -1}e^{-u/w\beta }/\Gamma(\alpha ),$ having
the mixing distribution ${G}_{2} ({\widehat{\beta}}_{1} ^{2};\lambda )$ with 
$\lambda $ = ${\beta}_{1} ^{2}/{\sigma }_{1}^{2} .$
\end{theorem}

\noindent\textbf{Proof.} Arguments follow step--by--step as in the proofs
given in Section 4.1. The details differ, but proceed similarly on noting
that $\bar{\boldsymbol{Y} }$ = $\QTR{up}{Diag}({n}_{1} ^{-1}\boldsymbol{1}_{{%
n}_{1} } ^{\prime},\ldots,{n}_{k} ^{-1}\boldsymbol{1}_{{n}_{k} } ^{\prime})%
\boldsymbol{Y} $ = $\boldsymbol{L}_{1} ^{\prime}\boldsymbol{Y} ,$ say,
whereas $\boldsymbol{R} $ = $\QTR{up}{Diag}(\boldsymbol{B}_{n_{1}} ,\ldots,%
\boldsymbol{B}_{n_{k}} )\boldsymbol{Y} $ = $\boldsymbol{L}_{2} ^{\prime}%
\boldsymbol{Y} .$ Conditional and unconditional moments follow directly as
linear functions, together with the idempotencies of $\{\boldsymbol{B}%
_{n_{i}} ;1 \leq i \leq k\}$ and the annihilations achieved through $\{%
\boldsymbol{B}_{n_{i}} \boldsymbol{1}_{{n}_{i} }  = \boldsymbol{0} ;1 \leq i
\leq k\}.$ A Gaussian model for $\mathcal{L}(\boldsymbol{Z} ),$ and Gaussian
errors during calibration, give conditional Gaussian laws for $\mathcal{L}(%
\bar{\boldsymbol{Y} }\!\mid\!{\widehat{\beta}}_{1} )$ and $\mathcal{L}(%
\boldsymbol{R} \!\mid\!{\widehat{\beta}}_{1} ),$ whereas unconditional
distributions are mixtures as in Section 2.2. Moreover, $\{{S}_{1}^{2}
,\ldots,{S}_{k}^{2} \}$ are conditionally independent given ${\widehat{\beta}%
}_{1} .$ As in the proof for Theorem 3(iv), the marginal density of $%
\mathcal{L}({\nu}_{i} {S}_{i}^{2} /{\omega }_{i}^{2} \!\mid\!{\widehat{\beta}%
}_{1} )$ is the scaled chi--squared density ${g}_{0} ({u}_{i} ;{\nu}_{i}
/2,2w)$ as defined preceding (2.2), with $w = {\widehat{\beta}}_{1} .$ Their
unconditional density now follows on mixing as in Section 2.2, as asserted
in conclusion (v), to complete our proof. $\square$

It is essential to examine effects of calibration in comparing variances
across the $k$ groups, typically based on $\{{S}_{1}^{2} ,\ldots,{S}_{k}^{2}
\}.$ In the ideal case where $({\beta}_{0} ,{\beta}_{1} )$ are known, we
would have $\{\QTR{up}{Var}({Y}_{ij} ) = {\beta}_{1} ^{2}{\omega }_{i}^{2}
;1 \leq i \leq k\},$ so that homoscedasticity across groups for measurements 
$\{{Y}_{ij} \}$ would be tantamount to that for $\{{Z}_{ij} \}.$ Under
calibration errors, however, Theorem 2(ii) shows that ${S}_{i}^{2} $
underestimates $\QTR{up}{Var}({Y}_{ij} ) = {\kappa }_{2} {\omega }_{i}^{2} 
+ {\sigma }_{0}^{2}  + {\sigma }_{1}^{2} {\mu}_{i} ^{2},$ the amount of
bias, ${B}_{i}  = -({\sigma }_{0}^{2}  + {\sigma }_{1}^{2} {\omega }_{i}^{2}
),$ being an artifact of the calibration process itself. Accordingly, it is
germane to examine homogeneity among the expected values $\{{\kappa }_{2} {%
\omega }_{1}^{2} ,\ldots,{\kappa }_{2} {\omega }_{k}^{2} \}$ of $\{{S}%
_{1}^{2} ,\ldots,{S}_{k}^{2} \}.$ To these ends let ${T}_{4} ({S}_{1}^{2}
,\ldots,{S}_{k}^{2} )$ be any scale--invariant statistic based on the sample
variances from the measurements $\{{Y}_{ij} ;1 \leq j \leq {n}_{i} ,1 \leq i
\leq k\}.$ The following result is basic.

\begin{theorem}
Let $\{{S}_{1}^{2} ,\ldots,{S}_{k}^{2} \}$ be within--sample variances from
the calibrated measurements $\{{Y}_{ij} ;1 \leq j \leq {n}_{i} ,1 \leq i
\leq k\}$ in a one--way experiment; let ${T}_{4} ({S}_{1}^{2} ,\ldots,{S}%
_{k}^{2} )$ be any scale--invariant statistic; and consider a Gaussian model
where $\{({Z}_{ij} -{\mu}_{i} )/{\omega }_{i} ;1 \leq j \leq {n}_{i} ,1 \leq
i \leq k\}$ are \textit{iid} ${N}_{1} (0,1)$ random variables independently
of $({\widehat{\beta}}_{0} ,{\widehat{\beta}}_{1} )$ under Gaussian
calibration errors. Then the distribution of ${T}_{4} ({S}_{1}^{2} ,\ldots,{S%
}_{k}^{2} )$ is identical to its normal--theory form, independently of $({%
\widehat{\beta}}_{0} ,{\widehat{\beta}}_{1} )$ and the empirical calibration
line.
\end{theorem}

\noindent\textbf{Proof.} The proof for Theorem 4(v) asserts $f(\boldsymbol{u}
;{\nu}_{1},\ldots,{\nu}_{k} ) = \prod_{i=1}^{k} {g}_{0} ({u}_{i} ;{\nu}_{i}
/2,2)$ as the conditional density for $\mathcal{L}({\nu}_{1} {S}_{1}^{2} /{%
\widehat{\beta}}_{1} {\omega }_{1}^{2} ,\ldots,{\nu}_{k} {S}_{k}^{2} /{%
\widehat{\beta}}_{1} {\omega }_{k}^{2} \!\mid\!{\widehat{\beta}}_{1} ),$
with ${\nu}_{i}  = {n}_{i} -1$ and ${g}_{0} (u;\alpha ,\beta )$ = $u^{\alpha
-1}e^{-u/\beta }/\beta ^{\alpha }\Gamma(\alpha ),$ so that $\{{S}_{1}^{2}
,\ldots,{S}_{k}^{2} \}$ are conditionally independent given ${\widehat{\beta}%
}_{1} .$ But since ${T}_{4} ({S}_{1}^{2} ,\ldots,{S}_{k}^{2} )$ is
scale--invariant, $\mathcal{L}[{T}_{4} ({S}_{1}^{2} ,\ldots,{S}_{k}^{2}
)\!\mid\!{\widehat{\beta}}_{1} ]$ = $\mathcal{L}[{T}_{4} ({S}_{1}^{2}
,\ldots,{S}_{k}^{2} )]$ unconditionally, to complete our proof. $\square$

It deserves note that meaningful comparisons among variances are necessarily
scale--invariant. Moreover, it is seen that procedures based on $\{{S}%
_{1}^{2} ,\ldots,{S}_{k}^{2} \}$ support tests for conditional hypotheses
that $\{\QTR{up}{Var}({Y}_{ij} \!\mid\!{\widehat{\beta}}_{1} ) = {\widehat{%
\beta}}_{1} ^{2}{\omega }^{2}  + {\sigma }_{0}^{2} ;1 \leq i \leq k\},$ or
equivalently, $H_{0}: {\omega }_{1}^{2}  = {\omega }_{2}^{2}  = \ldots = {%
\omega }_{k}^{2} ,$ to be tested against alternatives as appropriate.
Theorem 5 applies in the case of both null and non--null distributions of
invariant test statistics. Tests in common usage include

\begin{itemize}
\item Modifications of Bartlett's (1937) likelihood ratio test;

\item Cochran's (1941) test based on ${S}^{2} _{max}/({S}_{1}^{2}  + \ldots
+ {S}_{k}^{2} );$

\item Hartley's (1950) $F$--max test based on the maximal ratio $\QTR{up}{max%
}\{{S}_{i}^{2} /{S}_{j}^{2} \};$ and

\item Gnanadesikan's (1959) simultaneous comparisons of treatment variances
with a control.
\end{itemize}

To examine effects of calibration errors on the one--way analysis of
variance for comparing means, we proceed conditionally given ${\widehat{\beta%
}}_{1} ,$ first assuming that $\{\QTR{up}{Var}({Y}_{ij} \!\mid\!{\widehat{%
\beta}}_{1} ) = {\widehat{\beta}}_{1} ^{2}{\omega }^{2}  + {\sigma }_{0}^{2}
; 1 \leq j \leq {n}_{i} ,1 \leq i \leq k\},$ so that $V(\boldsymbol{Y}
\!\mid\!{\widehat{\beta}}_{1} )$ = ${\widehat{\beta}}_{1} ^{2}{\omega }^{2} 
\boldsymbol{I}_{n}  + {\sigma }_{0}^{2} \boldsymbol{1}_{n} \boldsymbol{1}%
_{n} ^{\prime}$ = $\boldsymbol{\Xi} ({\widehat{\beta}}_{1} )$ in the
notation of Section 3.2. We are concerned with comparative inferences
regarding elements of $\boldsymbol{\mu} ({\beta}_{1} )$ = ${\beta}_{1} [{\mu}%
_{1},\ldots,{\mu}_{k} ]^{\prime}$ from $E(\boldsymbol{Y} )$ = ${\beta}_{0} 
\boldsymbol{1}_{n}  + {\beta}_{1} [{\mu}_{1} \boldsymbol{1}_{{n}_{1} }
^{\prime},\ldots,{\mu}_{k} \boldsymbol{1}_{{n}_{k} } ^{\prime}]^{\prime}.$
Recall that $\boldsymbol{I}_{n} $ = $\boldsymbol{A}_{0}  + \boldsymbol{A}%
_{1}  + \boldsymbol{A}_{2} $ partitions $\boldsymbol{Y} ^{\prime}\boldsymbol{%
I}_{n} \boldsymbol{Y} $ = $\boldsymbol{Y} ^{\prime}\boldsymbol{A}_{0} 
\boldsymbol{Y}  + \boldsymbol{Y} ^{\prime}\boldsymbol{A}_{1} \boldsymbol{Y} 
+ \boldsymbol{Y} ^{\prime}\boldsymbol{A}_{2} \boldsymbol{Y} $ such that $%
\boldsymbol{Y} ^{\prime}\boldsymbol{A}_{0} \boldsymbol{Y} $ = $n\overline{Y}%
^{2},$ with $\overline{Y}$ as the grand mean and $\boldsymbol{A}_{0}  =
n^{-1}\boldsymbol{1}_{n} \boldsymbol{1}_{n} ^{\prime};$ $\boldsymbol{Y}
^{\prime}\boldsymbol{A}_{1} \boldsymbol{Y}  = \sum_{i=1}^{k}{n}_{i} ({%
\overline{Y}}_{i} -\overline{Y})^{2};$ and $\boldsymbol{Y} ^{\prime}%
\boldsymbol{A}_{2} \boldsymbol{Y} $ = $\sum_{i=1}^{k}\sum_{j=1}^{{n}_{i} }({Y%
}_{ij}  - {\overline{Y}}_{i} )^{2}.$ To validate the Fisher--Cochran theorem
conditionally requires that $\{\boldsymbol{A}_{i} \boldsymbol{\Xi} ({%
\widehat{\beta}}_{1} )\boldsymbol{A}_{j}  = \boldsymbol{0} ; i\neq j\}.$
Moreover, scale parameters associated with the quadratic forms are found as $%
\{{\kappa }_{i}^{2} \boldsymbol{A}_{i}  = \boldsymbol{A}_{i} \boldsymbol{\Xi}
({\widehat{\beta}}_{1} )\boldsymbol{A}_{i} ; i=1,2,3\},$ whereas
noncentrality parameters derive from expected mean squares. This program of
study is carried out next in support of the following.

\begin{theorem}
Let $\{{Y}_{ij}  = {\widehat{\beta}}_{0}  + {\widehat{\beta}}_{1} {Z}_{ij}
;1 \leq j \leq {n}_{i} ,1 \leq i \leq k\}$ be calibrated measurements in a
one--way experiment such that $\{({Z}_{ij} -{\mu}_{i} )/\omega ;1 \leq j
\leq {n}_{i} ,1 \leq i \leq k\}$ are \textit{iid} ${N}_{1} (0,1)$ random
variables independently of $({\widehat{\beta}}_{0} ,{\widehat{\beta}}_{1} )$
under Gaussian calibration errors.

(i) The analysis of variance test for equality of elements of $\boldsymbol{%
\mu} ({\beta}_{1} )$ = ${\beta}_{1} [{\mu}_{1},\ldots,{\mu}_{k} ]^{\prime},$
pertaining to the group measurement means, is identical in level and power
to its normal--theory form.

(ii) Supporting tests, based on linear contrasts among the group means, are
identical in level and power to their normal--theory forms.
\end{theorem}

\noindent\textbf{Proof.} To validate the Fisher--Cochran theorem
conditionally, observe $\{\boldsymbol{A}_{i} \boldsymbol{\Xi} ({\widehat{%
\beta}}_{1} )\boldsymbol{A}_{0}  = \boldsymbol{0} ; i=1,2\},$ since $%
\boldsymbol{A}_{i} \boldsymbol{\Xi} ({\widehat{\beta}}_{1} )\boldsymbol{A}%
_{0} $ = $\boldsymbol{A}_{i} ({\widehat{\beta}}_{1} ^{2}{\omega }^{2} 
\boldsymbol{I}_{n}  + {\sigma }_{0}^{2} \boldsymbol{1}_{n} \boldsymbol{1}%
_{n} ^{\prime})\boldsymbol{A}_{0} $ and $\{\boldsymbol{A}_{i} \boldsymbol{A}%
_{0}  = \boldsymbol{0} ; i=1,2\}.$ Similarly $\{\boldsymbol{A}_{i} 
\boldsymbol{\Xi} ({\widehat{\beta}}_{1} )\boldsymbol{A}_{j}  = \boldsymbol{0}
; (i,j)=1,2, i\neq j\},$ since $\boldsymbol{A}_{i} \boldsymbol{\Xi} ({%
\widehat{\beta}}_{1} )\boldsymbol{A}_{j} $ = $\boldsymbol{A}_{i} ({\widehat{%
\beta}}_{1} ^{2}{\omega }^{2} \boldsymbol{I}_{n}  + n{\sigma }_{0}^{2} 
\boldsymbol{A}_{0} )\boldsymbol{A}_{j} $ and $\{\boldsymbol{A}_{i} 
\boldsymbol{A}_{j}  = \boldsymbol{0} ; (i,j)=0,1,2, i\neq j\}$ from standard
properties of the one--way classification. Scale parameters, as determined
from $\{{\kappa }_{i}^{2} \boldsymbol{A}_{i}  = \boldsymbol{A}_{i} 
\boldsymbol{\Xi} ({\widehat{\beta}}_{1} )\boldsymbol{A}_{i} ; i=1,2\},$ are
found to be equal, namely $\{{\kappa }_{i}^{2} \boldsymbol{A}_{i}  = 
\boldsymbol{A}_{i} ({\widehat{\beta}}_{1} ^{2}{\omega }^{2} \boldsymbol{I}%
_{n}  + {\sigma }_{0}^{2} \boldsymbol{1}_{n} \boldsymbol{1}_{n} ^{\prime})%
\boldsymbol{A}_{i}  = {\widehat{\beta}}_{1} ^{2}{\omega }^{2} \boldsymbol{A}%
_{i} ; i=1,2\}$ from idempotency together with the annihilation $\{%
\boldsymbol{A}_{i} \boldsymbol{1}_{n}  = \boldsymbol{0} ; i=1,2\},$ so that $%
\{{\kappa }_{i}^{2}  = {\kappa }^{2} ; i=1,2\}.$ Finally, the noncentrality
parameters and degrees of freedom associated with $\{\boldsymbol{Y} ^{\prime}%
\boldsymbol{A}_{i} \boldsymbol{Y} ;i=1,2\}$ are determined from their
expected mean squares. These are $\{E(\boldsymbol{Y} ^{\prime}\boldsymbol{A}%
_{i} \boldsymbol{Y} \!\mid\!{\widehat{\beta}}_{1} ) = \QTR{up}{tr}(%
\boldsymbol{A}_{i} \boldsymbol{\Xi} ({\widehat{\beta}}_{1} ) + [\boldsymbol{%
\mu} ({\widehat{\beta}}_{1} )]^{\prime}\boldsymbol{A}_{i} \boldsymbol{\mu} ({%
\widehat{\beta}}_{1} ); i=1,2\}.$ It follows directly that 
\begin{equation}
E(\boldsymbol{Y} ^{\prime}\boldsymbol{A}_{1} \boldsymbol{Y} \!\mid\!{%
\widehat{\beta}}_{1} ) = \QTR{up}{tr}(\boldsymbol{A}_{1} \boldsymbol{\Xi} ({%
\widehat{\beta}}_{1} ) + [\boldsymbol{\mu} ({\widehat{\beta}}_{1} )]^{\prime}%
\boldsymbol{A}_{1} \boldsymbol{\mu} ({\widehat{\beta}}_{1} ) = (k-1){\kappa }%
^{2}  + {\widehat{\beta}}_{1} ^{2}\sum_{i=1}^{k}{n}_{i} ({\mu}_{i} -\bar{\mu}%
)^{2}  \notag
\end{equation}
\noindent with $\bar{\mu} = \sum_{i=1}^{k}{n}_{i} {\mu}_{i} /n.$ Similarly $%
E(\boldsymbol{Y} ^{\prime}\boldsymbol{A}_{2} \boldsymbol{Y} \!\mid\!{%
\widehat{\beta}}_{1} ) = \QTR{up}{tr}(\boldsymbol{A}_{2} \boldsymbol{\Xi} ({%
\widehat{\beta}}_{1} ) + [\boldsymbol{\mu} ({\widehat{\beta}}_{1} )]^{\prime}%
\boldsymbol{A}_{2} \boldsymbol{\mu} ({\widehat{\beta}}_{1} )$ = $(n-k){%
\kappa }^{2} $ since $[\boldsymbol{\mu} ({\widehat{\beta}}_{1} )]^{\prime}%
\boldsymbol{A}_{2} \boldsymbol{\mu} ({\widehat{\beta}}_{1} )$ = ${\widehat{%
\beta}}_{1} ^{2}\sum_{i=1}^{k}\sum_{j=1}^{{n}_{i} }({\mu}_{i} -{\mu}_{i}
)^{2} = 0.$ From these developments we infer that the distribution of the
ratio $F = (n-k)\boldsymbol{Y} ^{\prime}\boldsymbol{A}_{1} \boldsymbol{Y}
/(k-1)\boldsymbol{Y} ^{\prime}\boldsymbol{A}_{2} \boldsymbol{Y} $ satisfies $%
\mathcal{L}(F\!\mid\!{\widehat{\beta}}_{1} )$ = $F( k-1,n-k,\lambda ({%
\widehat{\beta}}_{1} ))$ with $\lambda ({\widehat{\beta}}_{1} )$ = ${%
\widehat{\beta}}_{1} ^{2}\sum_{i=1}^{k}{n}_{i} ({\mu}_{i} -\mu)^{2}/{%
\widehat{\beta}}_{1} ^{2}{\omega }^{2} $ = $\sum_{i=1}^{k}{n}_{i} ({\mu}_{i}
-\mu)^{2}/{\omega }^{2} .$ Thus the conditional and unconditional
distributions are identical, \textit{i.e.} $\mathcal{L}(F\!\mid\!{\widehat{%
\beta}}_{1} )$ = $\mathcal{L}(F)$ = $F(k-1,n-k,\lambda ),$ with $\lambda $ = 
$\sum_{i=1}^{k}{n}_{i} ({\mu}_{i} -\mu)^{2}/{\omega }^{2} .$ $\square$

\section{Diagnostics}

\label{S:dgnos}

\subsection{Objectives.}

\label{SS:objct} Calibration errors exact profound disturbances, both in
models and in data--analytic procedures, as shown. Myriad calibrated data
sets have been analyzed to date, supported of late by an evolving battery of
diagnostic tools. On these grounds it is tempting to dismiss the present
study as academic: For surely these issues long since would have surfaced in
practice, to be addressed accordingly. At issue is the capacity of known
diagnostics to uncover calibration--induced irregularities as documented
here. We now address these concerns with regard to induced correlations,
nonnormality, mixture distributions having excessive tails, and possible
outliers. For definiteness we return to the case of a single sample as in
Section 4.1.

\subsection{Correlation.}

\label{SS:detcor} Neither the conditional $({\sigma }_{0}^{2} /({\widehat{%
\beta}}_{1} ^{2}{\sigma }_{Z}^{2}  + {\sigma }_{0}^{2} ))$ nor the
unconditional $(({\sigma }_{0}^{2}  + {\sigma }_{1}^{2} {\mu}_{Z} ^{2})/({%
\kappa }_{2} {\sigma }_{Z}^{2}  + {\sigma }_{0}^{2}  + {\sigma }_{1}^{2} {\mu%
}_{Z} ^{2}))$ correlations need be negligible. Tests for correlation entail
dispersion matrices $V(\boldsymbol{Y} ) = \tau^{2}\boldsymbol{\Xi} (\omega ),
$ for which $\boldsymbol{\Xi} (\omega ) = (\boldsymbol{I}_{n}  + \omega 
\boldsymbol{A} )$ with $\boldsymbol{A} $ fixed and $\boldsymbol{\Xi} (\omega
)\in\mathbb{S}_{n}^{+}.$ Specializing gives $\tau^{2}\boldsymbol{\Xi}
(\omega )$ as $\boldsymbol{\Sigma } (\rho)$ under the equicorrelation models
encountered here. Tests of note are due to Durbin and Watson (1950, 1951,
1971), Anderson and Anderson (1950), Theil (1965), and others, all based on
versions of von Neumann's (1941) ratio $U = \boldsymbol{R} ^{\prime}%
\boldsymbol{B} \boldsymbol{R} /\boldsymbol{R} ^{\prime}\boldsymbol{R} ,$
with $\boldsymbol{R} $ as the observed residuals and with $\boldsymbol{B}
(n\times n)$ fixed. For further details see Kariya (1977). However, here the
unconditional distributions are all identical to their normal--theory forms
as if $\mathcal{L}(\boldsymbol{R} ) = {N}_{n} (\boldsymbol{0} ,{\sigma }^{2} 
\boldsymbol{B}_{n} ).$ This is seen from the proof for Theorem 3(ii), where $%
\mathcal{L}(\boldsymbol{R} \!\mid\!{\widehat{\beta}}_{1} )$ = ${N}_{n} (%
\boldsymbol{0} ,{\widehat{\beta}}_{1} ^{2}{\sigma }_{Z}^{2} \boldsymbol{B}%
_{n} ),$ together with the scale--invariance of $U = \boldsymbol{R} ^{\prime}%
\boldsymbol{B} \boldsymbol{R} /\boldsymbol{R} ^{\prime}\boldsymbol{R} ,$
assuring that $\mathcal{L}(\boldsymbol{R} ^{\prime}\boldsymbol{B} 
\boldsymbol{R} /\boldsymbol{R} ^{\prime}\boldsymbol{R} \!\mid\!{\widehat{%
\beta}}_{1} )$ = $\mathcal{L}(\boldsymbol{R} ^{\prime}\boldsymbol{B} 
\boldsymbol{R} /\boldsymbol{R} ^{\prime}\boldsymbol{R} )$ unconditionally.
All such diagnostics for correlative dependencies are totally blind, both to
the conditional $[V(\boldsymbol{Y} \!\mid\!{\widehat{\beta}}_{1} )$ = ${%
\widehat{\beta}}_{1} ^{2}{\sigma }_{Z}^{2} \boldsymbol{I}_{n}  + {\sigma }%
_{0}^{2} \boldsymbol{1}_{n} \boldsymbol{1}_{n} ^{\prime}]$ and unconditional 
$[V(\boldsymbol{Y} )$ = ${\kappa }_{2} {\sigma }_{Z}^{2} \boldsymbol{I}_{n} 
+ ({\sigma }_{0}^{2}  + {\sigma }_{1}^{2} {\mu}_{Z} ^{2})\boldsymbol{1}_{n} 
\boldsymbol{1}_{n} ^{\prime}]$ dispersion structures. In short, demonstrated
calibration--induced correlations cannot be discerned through conventional
diagnostic tools.

\subsection{Nonnormality.}

\label{SS:detnonnorm} Diagnostics for normality encompass both graphical and
hypothesis testing procedures. Graphics include plots of ordered residuals
against their normal--theory expectations. Common usage includes the scaled
residuals $\{{R}_{i} /{S}_{Y} ; 1 \leq i \leq n\},$ or the Studentized
residuals $\{{W}_{i} {R}_{i} /{S}_{Y} ; i=1, 2,\ldots,n\},$ standardized so
that $\QTR{up}{Var}({W}_{i} {R}_{i} ) = {\sigma }_{Y}^{2} .$ See Sections
2.12 and 5.7 of Myers (1990), for example. In calibrated data these residual
plots are indistinguishable from those for the conventional Gaussian model ${%
N}_{n} (\mu\boldsymbol{1}_{n} , {\sigma }^{2} \boldsymbol{I}_{n} ),$
whatever be the joint mixture density at (2.1) for the calibrated
measurements. This follows since $\mathcal{L}(\boldsymbol{R} /(\boldsymbol{R}
^{\prime}\boldsymbol{R} )^{\frac{1}{2}}\!\mid\!{\widehat{\beta}}_{1} )$ = $%
\mathcal{L}(\boldsymbol{R} /(\boldsymbol{R} ^{\prime}\boldsymbol{R} )^{\frac{%
1}{2}})$ from scale invariance, the latter as a scaled singular multivariate
Student's \textit{t}--distribution having $\nu = n-1$ degrees of freedom,
depending on neither ${\widehat{\beta}}_{1} $ nor ${\sigma }_{Y}^{2} .$

Tests for normality include the regression tests of Shapiro and Wilk (1965),
known to be powerful against a wide range of alternatives, to include skewed
or distributions having short or very long tails, even in small samples. See
Royston (1988), for example. These tests utilize statistics $W$ = $%
(\sum_{i=1}^{n}{w}_{i} {Y}_{[i]} )^{2}/(n-1){S}_{Y}^{2} ,$ where $\{{Y}%
_{[1]}  \leq {Y}_{[2]}  \leq \ldots \leq {Y}_{[n]} \}$ are the ordered
values of $\{{Y}_{1},\ldots,{Y}_{n} \},$ and $\{{w}_{1},\ldots,{w}_{n} \}$
are fixed weights. Such tests would appear promising for detecting the
nonstandard mixture distributions of calibrated measurements, where 
\begin{equation}
W = \frac{(\sum_{i=1}^{n}{w}_{i} {Y}_{[i]} )^{2}}{(n-1){S}_{Y}^{2} } = \frac{%
({\widehat{\beta}}_{0} \sum_{i=1}^{n}{w}_{i}  + {\widehat{\beta}}_{1}
^{2}\sum_{i=1}^{n}{w}_{i} {Z}_{[i]} )^{2}}{(n-1){S}_{Y}^{2} }.
\end{equation}
\noindent However, since $\sum_{i=1}^{n}{w}_{i}  = 0$ for these tests,
together with the identity ${S}_{Y}^{2}  = {\widehat{\beta}}_{1} ^{2}{S}%
_{Z}^{2} ,$ it follows that $W = (\sum_{i=1}^{n}{w}_{i} {Z}_{[i]} )^{2}/(n-1)%
{S}_{Z}^{2} .$ Then $\mathcal{L}(W\!\mid\!{\widehat{\beta}}_{1} ) = \mathcal{%
L}(W)$ holds unconditionally from cancellation. In short, all such
regression tests fail to distinguish between Gaussian distributions, and the
Gaussian mixtures of type (2.1). With regard to further variations on
regression tests, as in D'Agostino (1982), similar arguments show that none
is able to distinguish between Gaussian distributions and their mixtures
from calibrated measurements. Given sample moments $\{{m}_{r}  =
\sum_{i=1}^{n}({Y}_{i}  - \overline{Y})^{r}; r=2,3,4\},$ tests based on the
moment ratios $\{{b}_{1}  = {m}_{3}^{2} /{m}_{2}^{3} ,\,{b}_{2}  = {m}_{4} /{%
m}_{2}^{2} \}$ are useful against skewed alternatives or distributions
having excessive or short tails (D'Agostino (1982)). It is readily shown
that these ratios are precisely those obtainable from $\{{Z}_{1},\ldots,{Z}%
_{n} \},$ so that their null distributions are identical to those for which $%
\mathcal{L}(\boldsymbol{Y} ) = {N}_{n} (\mu\boldsymbol{1}_{n} , {\sigma }%
^{2} \boldsymbol{I}_{n} ),$ whatever be the joint mixture distribution as in
(2.1). On the other hand, the foregoing tests do offer a clear check on
normality of the distribution of $\{{Z}_{1},\ldots,{Z}_{n} \},$ on which the
mixtures (2.1) are predicated.

In short, conventional Gaussian diagnostics are bereft of any capacity to
distinguish between Gaussian errors, and Gaussian mixtures of type (2.1).
Thus radical calibration--induced departures from Gaussian models cannot be
discerned through routine screening using any of these diagnostics.

\subsection{Outliers.}

\label{SS:outdgnos} Commonly used diagnostics for a shift in location or
scale at observation ${Y}_{i} $ include the Studentized residuals ${t}_{i} 
= {R}_{i} /{S}_{Y} \sqrt{(1-1/n)},$ and the \textit{R--Student} deletion
diagnostic ${Rt}_{i}  = {R}_{i} /{S}_{-i} \sqrt{(1-1/n)},$ where ${S}_{-i} $
is the sample standard deviation found on deleting ${Y}_{i} $ from $\{{Y}%
_{1},\ldots,{Y}_{n} \}.$ As mixture distributions may have heavy tails, and
since conventional diagnostics for normality have failed, it is natural to
ask whether outlier diagnostics might be sensitive to observations from
mixtures of type (2.1). If so, then evidence for apparent outliers in
calibrated data instead might be attributable to the calibration process
itself. However, these diagnostics are all scale--invariant functions of the
observed residuals $\{{R}_{1},\ldots,{R}_{n} \},$ so that they are
indistinguishable from statistics derived from the standard Gaussian model ${%
N}_{n} (\mu\boldsymbol{1}_{n} ,{\sigma }^{2} \boldsymbol{I}_{n} ).$ In
short, conventional outlier diagnostics cannot distinguish between Gaussian
errors, and heavy--tailed mixtures as in (2.1), even if a shift in location
or scale has occurred at observation ${Y}_{i} .$

Section 5 has reexamined whether conventional diagnostics can detect
calibration--induced anomalies, to include correlations, nonnormality,
distributions having excessive tails, and possible outliers. Even radical
departures from conventional assumptions cannot be discerned through routine
screening using any of the aforementioned diagnostics. In summary, the
present study cannot be dismissed as merely academic, as evidence for
anomalies traceable to calibration could not have surfaced in practice
through a battery of diagnostic tools as it has evolved to date.

\section{Case Studies}

\label{S:CaseStud}

We apply the results of Section 4.1 to a numerical data set under the
assumptions of Theorem 3. Table 1 gives the percent of purity $(X)$ and the
octane number $(U)$ from a sample of $n=11$ different gasoline production
runs. Percent purity is determined readily, whereas octane numbers require
expensive and time--consuming dynamic laboratory tests; hence the need for
calibration. \newpage

\begin{center}
\textbf{Table 1.} Percent of purity $(X)$ and octane number $(U)$ of gasoline

\vspace{0.1in} 
\begin{tabular}{|c|c|c|c|c|c|c|c|c|c|c|c|}
\hline
\hspace{0.1in}$X$\hspace{0.1in} & 99.8 & 99.7 & 99.6 & 99.5 & 99.4 & 99.3 & 
99.2 & 99.1 & 99.0 & 98.9 & 98.8 \\ \hline
$U$ & 88.6 & 86.4 & 87.2 & 88.4 & 87.2 & 86.8 & 86.1 & 87.3 & 86.4 & 86.6 & 
87.1 \\ \hline
\end{tabular}%
%
\end{center}

\noindent The least--squares fit for $U={\beta}_{0}  + {\beta}_{1} (X-%
\overline{X}) + \varepsilon $ has $\{n=11,{\widehat{\beta}}_{0} =87.2818,{%
\widehat{\sigma }}_{0} =0.1846,{\widehat{\beta}}_{1} =1.8546, {\widehat{%
\sigma }}_{1} =0.5837\}.$ Suppose that subsequent determinations of percent
purity satisfy $\{\mathcal{L}({Z}_{i} ) = {N}_{1} (0, 1); 1 \leq i \leq n\},$
so that calibrated measurements are recovered as $\{{Y}_{i}  = {\widehat{%
\beta}}_{0}  + {\widehat{\beta}}_{1} {Z}_{i} ;1 \leq i \leq n\}$ in units of
octane number. Then the distribution of $\overline{Y}$ is the mixture of a
normal distribution ${N}_{1} (\mu(t),\Xi(t)),$ with $\mu(t) = {\beta}_{0}  +
t{\mu}_{Z} $ and $\Xi(t) = {t}^{2} {\sigma }_{Z}^{2} /n+{\sigma }_{0}^{2} ,$
having the density ${g}_{1} (u;\mu(t),\Xi(t)),$ with mixing distribution ${N}%
_{1} ({\beta}_{1} ,{\sigma }_{1}^{2} )$ having the density $d{G}_{1} (t).$\
For convenience, we write this as $\mathcal{L}(\overline{Y}) = {N}_{1} ({%
\beta}_{0}  + t{\mu}_{Z} ,{t}^{2} {\sigma }_{Z}^{2} /n+{\sigma }_{0}^{2} ){%
\Lambda }_{t} {N}_{1} ({\beta}_{1} ,{\sigma }_{1}^{2} ),$ where ${\Lambda }%
_{t} $ designates the mixing operation. Accordingly, the density of $%
\overline{Y}$ is 
\begin{align}
{f}_{1} (u) =&\int\limits_{-\infty }^{\infty }{g}_{1} (u;{\beta}_{0}  + t{\mu%
}_{Z} ,{t}^{2} {\sigma }_{Z}^{2} /n+{\sigma }_{0}^{2} )d{G}_{1} (t)  \notag
\label{density for mean(Y)} \\
=&\int\limits_{-\infty }^{\infty }\frac{\exp \left( -\frac{[u-({\beta}_{0} 
+ t{\mu}_{Z} )]^{2}}{2({t}^{2} {\sigma }_{Z}^{2} /n+{\sigma }_{0}^{2} )}%
\right) }{\sqrt{2\pi ({t}^{2} {\sigma }_{Z}^{2} /n+{\sigma }_{0}^{2} )}}%
\frac{\exp \left( - \frac{(t-{\beta}_{1} )^{2}}{2{\sigma }_{1}^{2} }\right) 
}{\sqrt{2\pi {\sigma }_{1}^{2} }}dt
\end{align}%
a function of the parameters $\Omega$ = $\{n,{\beta}_{0} ,{\sigma }_{0} ,{\mu%
}_{Z} ,{\sigma }_{Z} ,{\beta}_{1} ,{\sigma }_{1} \},$ with skewness $%
0.1464\times 10^{-6}$ and kurtosis $3.855,$ and with conditional mean $E(%
\overline{Y}) = 87.2818,$ given the empirical calibration.

Using equation (\ref{density for mean(Y)}), we compute the $95\%$
probability region for $\overline{Y}$ as $(86.037,88.526)$ compared to ${%
\widehat{\beta}}_{0}  \pm 1.96|{\widehat{\beta}}_{1} |{\sigma }_{Z} /\sqrt{11%
} = (86.184,88.376)$ if $\overline{Y}$ were normal. This latter interval is
actually a $92.2\%$ probability region. In addition, the density ${f}_{1} (u)
$ is bell-shaped but is not normal. The Table 2 gives its moments (mean,
variance), and moment ratios (skewness $(\gamma ),$ kurtosis $(\kappa )$)
for selected values of the parameters $\Omega $.

The scaled sample variance is a mixture of a gamma distribution, ${G}_{0}
(\cdot, \cdot),$ with mixing distribution $d{G}_{2} (w)$ as a non-central
chi-squared distribution, to give $\mathcal{L}((n-1){S}_{Y}^{2} /{\sigma }%
_{1}^{2} {\sigma }_{Z}^{2} ) = {G}_{0} ((n-1)/2,\,2t){\Lambda }_{t} {\chi}%
_{1}^{2} (\lambda  = {\beta}_{1}^{2} /{\sigma }_{1}^{2} )$ with $E({S}%
_{Y}^{2} /{\sigma }_{Z}^{2} {\sigma }_{1}^{2} ) = (1 + \lambda )$ which for
the Octane Data is the conditional value $E({S}_{Y}^{2} )/{\sigma }_{Z}^{2} {%
\sigma }_{1}^{2}  = 11.095,$ so that $E({S}_{Y}^{2} ) = 3.780.$ The mixture
distribution has density 
\begin{align}
{f}_{0} (u) =&\frac{u^{\nu/2-1}}{2^{\nu/2}\Gamma (\nu/2)}\int\limits_{0}^{%
\infty}w^{-\nu/2}e^{-u/2w}d{G}_{2} (w)  \notag  \label{density for S_sq} \\
=&\frac{u^{\nu/2-1}}{2^{\nu/2}\Gamma (\nu/2)}\int\limits_{0}^{\infty}
w^{-\nu/2}e^{-u/2w}[\frac{e^{-\lambda /2 -u/2}}{2^{1/2}} \sum\limits_{j=0}^{%
\infty}(\frac{\lambda }{4})^{j}\frac{u^{j-1/2}}{j!\Gamma (1/2+j)}]dw
\end{align}%
with $\nu =(n-1).$

If $Y$ were normal with ${\widehat{\beta}}_{1} $ a constant, then a 95\%
probability region for ${S}_{Y}^{2} $ could be found from $P[{\chi}^{2}
(10;0.025) < (n-1){S}_{Y}^{2} /{\widehat{\beta}}_{1} ^{2}{\sigma }_{Z}^{2} 
< {\chi}^{2} (10;0.975)]$ or equivalently $P[1.1167 < {S}_{Y}^{2}  <7.0449]$
which actually is a $74\%$ probability region when variation in ${\widehat{%
\beta}}_{1} $ is taken into account. The correct probability region is found
by numerically integrating Equation (\ref{density for S_sq}) to get $P[10.8
< (n-1){S}_{Y}^{2} /{\sigma }_{1}^{2} {\sigma }_{Z}^{2}  < 336.5]$ = 0.95 = $%
P[0.3680 < {S}_{Y}^{2}  < 11.46].$ We find that using the first 20 terms in
the infinite sum is adequate.


\begin{center}
\textbf{Table 2.} The moments (mean, variance), and moment ratios (skewness $%
(\gamma ),$ kurtosis $(\kappa )$) of $\mathcal{L}(\overline{Y})$ for
selected values of the parameters $\Omega$ = $\{n,{\beta}_{0} ,{\sigma }_{0}
,{\mu}_{Z} ,{\sigma }_{Z} , {\beta}_{1} ,{\sigma }_{1} \}$

\vspace{0.1in} 
\begin{tabular}{|l|l|l|l|l|l|l|l|l|l|l|l|}
\hline
$n$ & ${\beta}_{0} $ & ${\sigma }_{0} $ & ${\mu}_{Z} $ & ${\sigma }_{Z} $ & $%
{\beta}_{1} $ & ${\sigma }_{1} $ &  & $E(\overline{Y})$ & $Var(\overline{Y})$
& $\gamma  $ & $\kappa  $ \\ \hline
10 & 1 & 1 & 1 & 1 & 1 & 1 &  & 2.0000 & 2.2000 & 0.1839 & 3.2851 \\ \hline
20 & 1 & 1 & 1 & 1 & 1 & 1 &  & 2.0000 & 2.1000 & 0.0986 & 3.1463 \\ \hline
&  &  &  &  &  &  &  &  &  &  &  \\ \hline
20 & .5 & 1 & 1 & 1 & 1 & 1 &  & 1.5000 & 2.1000 & 0.0986 & 3.1463 \\ \hline
20 & 2 & 1 & 1 & 1 & 1 & 1 &  & 3.0000 & 2.1000 & 0.0986 & 3.1463 \\ \hline
&  &  &  &  &  &  &  &  &  &  &  \\ \hline
20 & 1 & .5 & 1 & 1 & 1 & 1 &  & 2.0000 & 1.3500 & 0.1913 & 3.3539 \\ \hline
20 & 1 & 2 & 1 & 1 & 1 & 1 &  & 2.0000 & 5.100 & 0.0260 & 3.0248 \\ \hline
&  &  &  &  &  &  &  &  &  &  &  \\ \hline
20 & 1 & 1 & .5 & 1 & 1 & 1 &  & 1.5000 & 1.3500 & 0.0956 & 3.1070 \\ \hline
20 & 1 & 1 & 2 & 1 & 1 & 1 &  & 3.0000 & 5.1000 & 0.0521 & 3.0940 \\ \hline
&  &  &  &  &  &  &  &  &  &  &  \\ \hline
20 & 1 & 1 & 1 & .5 & 1 & 1 &  & 2.0000 & 2.0250 & 0.0260 & 3.0373 \\ \hline
20 & 1 & 1 & 1 & 2 & 1 & 1 &  & 2.0000 & 2.4000 & 0.3327 & 3.5417 \\ \hline
&  &  &  &  &  &  &  &  &  &  &  \\ \hline
20 & 1 & 1 & 1 & 1 & .5 & 1 &  & 1.5000 & 2.0625 & 0.0506 & 3.1463 \\ \hline
20 & 1 & 1 & 1 & 1 & 2 & 1 &  & 3.0000 & 2.2500 & 0.1778 & 3.1452 \\ \hline
&  &  &  &  &  &  &  &  &  &  &  \\ \hline
20 & 1 & 1 & 1 & 1 & 1 & .5 &  & 2.0000 & 1.3125 & 0.0499 & 3.0267 \\ \hline
20 & 1 & 1 & 1 & 1 & 1 & 2 &  & 2.0000 & 5.2500 & 0.0998 & 3.3614 \\ \hline
&  &  &  &  &  &  &  &  &  &  &  \\ \hline
10 & 1 & .5 & 1 & 2 & 1 & 2 &  & 2.0000 & 6.2500 & 0.6144 & 5.5559 \\ \hline
\end{tabular}%
%
\end{center}

\vspace{0.1in}

For the density of ${t}_{0}^{2}  = n(\overline{Y}-{\mu}_{Y}^{0} )^{2}/{S}%
_{Y}^{2} ,$ set 
\begin{equation}
\{\nu = n-1,\,\delta = ({\mu}_{Y} -{\mu}_{Y}^{0} )^{2}/{\sigma }_{1}^{2} {%
\sigma }_{Z}^{2} ,\,\delta({\widehat{\beta}}_{1} ) = \delta/({\widehat{\beta}%
}_{1} ^{2}/{\sigma }_{1}^{2} )\}.  \notag
\end{equation}
Its density is found on mixing the non-central ${t}^{2} (\nu ,\delta /t)$
over a non-central chi-squared as the mixing distribution, which we write as 
$\mathcal{L}({t}_{0}^{2} )={t}^{2} (\nu ,\delta /t){\Lambda }_{t} {\chi}%
_{1}^{2} (\lambda  ={\beta}_{1}^{2} /{\sigma }_{1}^{2} ).$ The density is
given by 
\begin{align}
&g(u;\nu ,\delta ,\lambda ) =  \notag \\
&\int_{0}^{\infty} \frac{1}{\nu }\sum_{j=0}^{\infty }\frac{(\frac{ \delta /t%
}{2})^{j}e^{-\frac{\delta /t}{2}}(\frac{u}{\nu })^{-\frac{1}{2}+j}}{j!B(%
\frac{1+2j}{2},\frac{\nu}{2})(1+\frac{u}{\nu })^{\frac{1}{2}+\frac{\nu }{2}%
+j}} \frac{e^{-\frac{\lambda  }{2}-\frac{t}{2}}}{2^{\frac{1}{2}}}%
\sum_{k=0}^{\infty }\frac{(\frac{\lambda  }{4})^{k}t^{k-\frac{1}{2}}}{%
k!\Gamma (\frac{1+2k}{2})}dt.
\end{align}%
For the first sum, we use the first ${N}_{1} =15$ terms, and for the second
sum the first ${N}_{2} =30$ terms.

This distribution is useful for computing the power of the ${t}^{2} $ test.
For example, with $n=11$ so $\nu =10$, the $95\%$ critical value is 4.9646
with $\delta =0.$ Table 3 gives the power of the test for $\delta
=\{0,1,4,9\}$ and $\lambda  = \{1,4,9\}.$

\vspace{0.1in}

\begin{center}
\textbf{Table 3.} Power for the test $H_{0}:{\mu}_{Y}  ={\mu}_{Y}^{0} $
against $H_{1}:{\mu}_{Y}  \neq {\mu}_{Y}^{0} $\\[0pt]
for $\delta\in\{0,1,4,9\}$ and $\lambda \in\{1,4,9\}$

\vspace{0.1in} 
\begin{tabular}{|c|c|c|c|}
\hline
\hspace{0.2in}$\lambda  =$\hspace{0.3in} & 1 & 4 & 9 \\ \hline
\hspace{0.2in}$\delta =0$\hspace{0.2in} & .950 & 950 & 950 \\ \hline
\hspace{0.2in}$\delta =1$\hspace{0.2in} & .691 & .863 & .928 \\ \hline
\hspace{0.2in}$\delta =4$\hspace{0.2in} & .485 & .742 & .876 \\ \hline
\hspace{0.2in}$\delta =9$\hspace{0.2in} & .329 & .608 & .799 \\ \hline
\end{tabular}%
%
\end{center}

\vspace{0.1in}

For the Octane Data, the power of the test of $H_{0}:{\mu}_{Y}  = {\mu}%
_{Y}^{0} ,$ with $({\mu}_{Y}  - {\mu}_{Y}^{0} )^{2}=1,$ has $\nu =10,$ $%
\delta =({\mu}_{Y}  - {\mu}_{Y}^{0} )^{2}/{\sigma }_{Z}^{2} {\sigma }%
_{1}^{2} $ = $[(1)(0.5837)]^{-2} = 2.9351,$ and $\lambda $ = ${\beta}%
_{1}^{2} /{\sigma }_{1}^{2} $= $(1.8546/0.5837)^{2} = 10.0953.$ The power of
the test is 90$\%.$

An equivalent form for $\mathcal{L}({t}_{0}^{2} )$ is based on ${t}_{0} $ as 
$\mathcal{L}({t}_{0} )=t(\nu,{\delta }_{0} /s){\Lambda }_{s} \sqrt{{\chi}%
_{1}^{2} (\lambda  = {\beta}_{1}^{2} /{\sigma }_{1}^{2} )},$ mixing over a
shifted half--normal distribution. Its density is%
\begin{align}
&f_{{t}_{0} }(u) =  \notag \\
&\int_{0}^{\infty}\frac{e^{-\frac{({\delta }_{0} /s)^{2}}{2}}\Gamma (( \nu +
1)/2)(\frac{\nu}{\nu +u^{2}})^{\frac{\nu}{2}+\frac{1}{2}}}{\sqrt{\pi\nu}%
\Gamma (\nu/2)}\sum_{j=0}^{\infty}\frac{\Gamma ((\nu + j + 1)/2)}{j!\Gamma
((\nu + 1)/2)}\left(\frac{\sqrt{2}u{\delta }_{0} /s}{\sqrt{\nu +u^{2}}}%
\right)^{j}d\sqrt{{G}_{2} }(s)
\end{align}%
having non-centrality parameter ${\delta }_{0} =\sqrt{\delta}$ and 
\begin{equation*}
d\sqrt{{G}_{2} }(s)=\frac{e^{-(s-{\lambda }_{0} )^{2}/2}+e^{-(-s-{\lambda }%
_{0} )^{2}/2}}{\sqrt{2\pi }}ds 
\end{equation*}%
with non-centrality parameter ${\lambda }_{0} =\sqrt{\lambda }.$ This series
has faster convergence and we used $N=20$ terms in the forgoing power
calculations with noncentrality parameters ${\delta }_{0}  = \sqrt{\delta}$
and ${\lambda }_{0}  = \sqrt{\lambda }$. Computations reported here were
executed by the second author using the Maple software package.

\section{Conclusions}

\label{S:concl}

In summary, the widespread and necessary use of calibration may have
devastating effects, even on elementary data--analytic procedures pertaining
to location and scale parameters. It is unfortunate that these difficulties
cannot be flagged by the ever expanding use of available diagnostic tools.
It thus is incumbent on knowledgeable users of statistical methodology, and
the statistical consultants advising them, to assess the extent of these
difficulties as they might impact the analysis and interpretation of data in
a particular experimental setting. Let the user be forewarned. Fortunately,
comparisons among means and among variances, in the context of comparative
one--way experiments, are largely unaffected by the use of calibrated
instruments when subject to errors of calibration, provided that the results
are interpreted accordingly.

\section{APPENDIX}

\label{S:appendx}

It is germane to examine the comparative concentration of probability
measures on $\mathbb{R}^{n}.$ Following Sherman (1955), the measure $%
\mu(\cdot)$ is said to be \textit{more peaked} about $\boldsymbol{0} \in%
\mathbb{R}^{n}$ than $\nu(\cdot)$ if and only if $\mu(A)\ge\nu(A)$ for every
set $A$ in the class $\mathcal{C}(n)$ comprising the convex sets in $\mathbb{%
R}^{n}$ symmetric under reflection through $\boldsymbol{0} \in\mathbb{R}^{n}.
$ For scale mixtures of Gaussian measures on $\mathbb{R}^{n},$ their
peakedness ordering is tantamount to the stochastic ordering of their mixing
distributions. Details follows.

\begin{lemma}
Let ${GM}_{n} (\boldsymbol{\theta} ,\boldsymbol{\Xi} ,{G}_{1} )$ and ${GM}%
_{n} (\boldsymbol{\theta} ,\boldsymbol{\Xi} ,{G}_{2} )$ be Gaussian mixtures
on $\mathbb{R}^{n}$ of type (2.2) having mixing distributions ${G}_{1}
(\cdot)$ and ${G}_{2} (\cdot)$ on $\mathbb{R}^{1}_{+}.$ Then ${GM}_{n} (%
\boldsymbol{\theta} ,\boldsymbol{\Xi} ,{G}_{1} )$ is more peaked about $%
\boldsymbol{\theta} \in\mathbb{R}^{n}$ than ${GM}_{n} (\boldsymbol{\theta} ,%
\boldsymbol{\Xi} ,{G}_{2} )$ if and only if ${G}_{1} (t)\leq{G}_{2} (t)$ for
every $t>0.$
\end{lemma}

\noindent\textbf{Proof.} The ordering ${G}_{1} (t)\leq{G}_{2} (t),$ \textit{%
i.e.,} that ${G}_{1} (\cdot)$ is stochastically larger than ${G}_{2} (\cdot),
$ holds if and only if there are increasing functions $\{{\psi}_{1} (\cdot),{%
\psi}_{2} (\cdot)\},$ ordered pointwise as ${\psi}_{1} (t)\geq{\psi}_{2} (t),
$ together with a random variable $U,$ such that ${G}_{1} (t) = P({\psi}_{1}
(U)\leq t)$ and ${G}_{2} (t)$ = $P({\psi}_{2} (U)\leq t);$ see Lemma 1, page
84 of Lehmann (1986), for example. Accordingly, we provisionally write $%
\mu(A)$ = $\int_{A}f(\boldsymbol{x} ;\boldsymbol{\theta} ,\boldsymbol{\Xi} ,{%
G}_{1} )d\boldsymbol{x} $ and $\nu(A)$ = $\int_{A}f(\boldsymbol{x} ;%
\boldsymbol{\theta} ,\boldsymbol{\Xi} ,{G}_{2} )d\boldsymbol{x} ,$ and their
difference as 
\begin{equation}
\mu(A) - \nu(A) = \int_{0}^{\infty}\int_{A}[g(\boldsymbol{x} ;\boldsymbol{%
\theta} ,\boldsymbol{\Xi} /{\psi}_{1} (t)) - g(\boldsymbol{x} ;\boldsymbol{%
\theta} ,\boldsymbol{\Xi} /{\psi}_{2} (t))]d\boldsymbol{x}  dG(t).  \notag
\end{equation}
\noindent Given that ${G}_{1} (t)\leq{G}_{2} (t),$ so that ${\psi}_{1}
(t)\geq{\psi}_{2} (t),$ the ordering $\int_{A}[g(\boldsymbol{x} ;\boldsymbol{%
\mu} ,\boldsymbol{\Xi} /{\psi}_{1} (t)) - g(\boldsymbol{x} ;\boldsymbol{\mu}
,\boldsymbol{\Xi} /{\psi}_{2} (t))]d\boldsymbol{x}  \geq 0$ follows
point--wise for each fixed $t\in\mathbb{R}^{1}_{+}$ from Corollary 3 of
Anderson (1955), since $\boldsymbol{\Xi} /{\psi}_{2} (t)\succeq_{L}%
\boldsymbol{\Xi} /{\psi}_{1} (t))$ uniformly in $t.$ That $%
[\mu(A)-\nu(A)]\geq 0$ now follows directly. Conversely, suppose that $%
\mu(A)\geq\nu(A).$ We now apply the converse to Anderson's (1955) Corollary
3, as proved in Jensen (1984), to infer that ${\psi}_{1} (t)\geq{\psi}_{2}
(t)$ for each $t>0,$ thus establishing the necessity of the condition ${G}%
_{1} (t)\leq{G}_{2} (t),$ to complete our proof. $\square$

\end{document}